\documentclass[10pt]{amsart}
\usepackage{latexsym}   
\usepackage{amssymb}    
\usepackage{amsmath}    
\usepackage{amsthm}     
\usepackage[all]{xy}
\usepackage{hyperref}
\newcommand{\pa}{\partial}\newcommand{\al}{\alpha}
\newcommand{\be}{\beta}
\newcommand{\ga}{\gamma}
\newcommand{\Ga}{\Gamma}\newcommand{\del}{\delta}
\newcommand{\Del}{\Delta}

\newcommand{\la}{\lambda}\newcommand{\La}{\Lambda}\newcommand{\om}{\omega}

\newcommand{\ti}{\tilde}

\renewcommand{\thefootnote}

\title[ Bianchi's B\"{a}cklund transformation for higher dimensional quadrics]
{ Bianchi's B\"{a}cklund transformation for higher dimensional
quadrics}

\author[  Ion I. Dinc\u{a}]{Ion I. Dinc\u{a}}
\address{Faculty of Mathematics and Informatics,
University of Bucharest,  14 Academiei Str., 010014, Bucharest,
Romania}
 \email{dinca@gta.math.unibuc.ro}
\thanks{Supported by the University of Bucharest}

\begin{document}

\keywords{B\"{a}cklund transformation, Bianchi Permutability
Theorem, (confocal) quadrics, common conjugate systems, (discrete)
deformations in $\mathbb{C}^{2n-1}$ of quadrics in
$\mathbb{C}^{n+1}$, {\it The Method} of Archimedes}

\begin{abstract}
We provide a generalization of Bianchi's B\"{a}cklund
transformation from $2$-dimensional quadrics to higher dimensional
quadrics. The starting point of our investigation is the higher
dimensional (infinitesimal) version of Bianchi's main four
theorems on the theory of deformations of quadrics and Bianchi's
treatment of the B\"{a}cklund transformation for diagonal
paraboloids via conjugate systems.
\end{abstract}

\maketitle

\tableofcontents \pagenumbering{arabic}

\section{Introduction}

In 1859 the French Academy posed the problem:
\newline
\begin{center}
{\it To find all surfaces applicable to a given one}.
\newline
\end{center}
It became the driving force which led to the flourishing of the
classical differential geometry in the second half of the
XIX$^{\mathrm{th}}$ century and its profound study by illustrious
geometers led to interesting results (see Bianchi
\cite{B1},\cite{B2},\cite{B3}, Darboux \cite{D1}, Eisenhart
\cite{E1},\cite{E2}, Sabitov \cite{S1} and its references for
results up to 1990's or Spivak (\cite{S2}, Vol {\bf 5})). Today it
is still an open problem in its full generality, but basic
familiar results like the Gau\ss-Bonnet Theorem and the
Codazzi-Mainardi equations (independently discovered also by
Peterson) were first communicated to the French Academy. A list
(most likely incomplete) of the winners of the prize includes
Bianchi, Bonnet, Guichard, Weingarten.

Up to 1899 deformations of the (pseudo-)sphere and isotropic
quadrics without center (from a metric point of view they can be
considered as metrically degenerate quadrics without center)
together with their {\it B\"{a}cklund} (B) transformation and the
complementary transformation of deformations of surfaces of
revolution were investigated by geometers such as B\"{a}cklund,
Bianchi, Darboux, Goursat, Hazzidakis, Lie, Weingarten, etc.

In 1899 Guichard discovered that when quadrics with(out) center
and of revolution around the focal axis roll on their deformations
their foci describe constant mean curvature (minimal) surfaces
(and Bianchi proved the converse: all constant mean curvature
(minimal) surfaces can be realized in this way).

With Guichard's result the race to find the deformations of
general quadrics was on; it ended with Bianchi's discovery
\cite{B0} from 1906 of the B transformation for quadrics and the
{\it applicability correspondence provided by the Ivory affinity}
(ACPIA).

Note also that Peterson's work on deformations of general quadrics
preceded that of Bianchi, Calapso, Darboux, Guichard and
\c{T}i\c{t}eica's from the years 1899-1906 by two decades, but
unfortunately most of his works (including his independent
discovery of the Codazzi-Mainardi equations and of the
Gau\ss-Bonnet Theorem) were made known to Western Europe mainly
after they were translated in 1905 from Russian to French (as is
the case with his deformations of quadrics \cite{P}, originally
published in 1883 in Russian).

The work of these illustrious geometers on deformations in
$\mathbb{C}^3$ of quadrics in $\mathbb{C}^3$ (there is no other
class of surfaces for which an interesting theory of deformation
has been built) is one of the crowning achievements of the golden
age of classical geometry of surfaces (Darboux split the prize of
the French Academy between Bianchi and Guichard for solving the
problem for quadrics (Guichard had other results including the G
transformation, but the G transformation turned out to be the
composition of two B transformations); it was solved in the sense
that solutions depending on arbitrarily many constants are
produced by algebraic procedures once an $1$-dimensional family of
Riccati equations is presumed solved) and at the same time it
opened new areas of research (such as affine and projective
differential geometry) continued later by other illustrious
geometers (Blaschke, Cartan, etc).

Note that Calapso in \cite{Ca} has put Bianchi's B transformation
of deformations in $\mathbb{C}^3$ of $2$-dimensional quadrics with
center in terms of common conjugate systems (the condition that
the conjugate system on a $2$-dimensional quadric is a conjugate
system on one of its deformations in $\mathbb{C}^3$ was known to
Calapso for a decade, but the B transformation for quadrics eluded
Calapso since the common conjugate system was best suited for this
transformation only at the analytic level: to find the
applicability correspondence one must renounce the correspondence
of the common conjugate systems provided by the {\it Weingarten}
(W) congruence in favor of the ACPIA).

In 1919-20 Cartan has shown in \cite{C} (using mostly projective
arguments and his exterior differential systems in involution and
exteriorly orthogonal forms tools) that space forms of dimension
$n$ admit rich families of deformations (depending on $n(n-1)$
functions of one variable) in surrounding space forms of dimension
$2n-1$, that such deformations admit lines of curvature (given by
a canonical form of exteriorly orthogonal forms; thus they have
flat normal bundle; since the lines of curvature on
$n$-dimensional space forms (when they are considered by
definition as quadrics in surrounding $(n+1)$-dimensional space
forms) are undetermined, the lines of curvature on the deformation
and their corresponding curves on the quadric provide the common
conjugate system) and that the codimension $n-1$ cannot be lowered
without obtaining rigidity as the deformation being the defining
quadric.

In 1979, upon a suggestion from S. S. Chern and using Chebyshev
coordinates on $\mathbb{H}^n(\mathbb{R})$ (the Cartan-Moore
Theorem; they are the lines of curvature on and thus in bijective
correspondence with deformations of $\mathbb{H}^n(\mathbb{R})$ in
$\mathbb{R}^{2n-1}$) Tenenblat-Terng have developed in \cite{TT1}
the B transformation of $\mathbb{H}^n(\mathbb{R})$ in
$\mathbb{R}^{2n-1}$ (and Terng in \cite{T2} has developed the {\it
Bianchi Permutability Theorem} (BPT) for this B transformation).

In 1983 Berger, Bryant and Griffiths \cite{B} proved (including by
use of tools from algebraic geometry) in particular that Cartan's
essentially projective arguments (including the exterior part of
his exteriorly orthogonal forms tool) can be used to generalize
his results to $n$-dimensional general quadrics with positive
definite linear element (thus they can appear as quadrics in
$\mathbb{R}^{n+1}$ or as space-like quadrics in
$\mathbb{R}^n\times(i\mathbb{R})$) admitting rich families of
deformations (depending on $n(n-1)$ functions of one variable) in
surrounding Euclidean space $\mathbb{R}^{2n-1}$, that the
codimension $n-1$ cannot be lowered without obtaining rigidity as
the deformation being the defining quadric and that quadrics are
the only Riemannian $n$-dimensional manifolds that admit a family
of deformations in $\mathbb{R}^{2n-1}$ as rich as possible for
which the exteriorly orthogonal forms tool (naturally appearing
from the Gau\ss\ equations) can be applied.

The starting point of our investigation is the higher dimensional
(infinitesimal) version of Bianchi's main four theorems on the
theory of deformations of quadrics and Bianchi's treatment of the
B transformation for paraboloids via conjugate systems.

Similarly to Bianchi's original (pre-Ivory affinity) approach to
the deformation problem for quadrics, when he made a link between
deformations of diagonal paraboloids and the sine-Gordon equation,
we shall consider first paraboloids, since the higher dimensional
version of the sine-Gordon equation (namely Terng's {\it
generalized sine-Gordon equation} (GSGE)) together with its B
transformation is already completed by Tenenblat-Terng and Terng
has already completed the BPT in this case.

Again just like Bianchi considered {\it isothermic-conjugate}
coordinates on the quadric (that is the second fundamental form is
a multiple of the identity) coupled with the conjugate system
common to the quadric and its deformation (the change between
these conjugate systems provides the angle that satisfies the
sine-Gordon equation and the ACPIA is best seen at the level of
the initial isothermic-conjugate system on the quadric), we shall
consider isothermic-coordinates on the quadric coupled with the
conjugate system common to both the quadric and its deformation,
since this interplay between the two systems of coordinates will
make the deformation problem amenable to an attack strategy (the
common conjugate system property provides the canonical form of
exteriorly orthogonal forms). The B transformation will turn out
to be similar in some aspects to that of Tenenblat-Terng (in fact
at the level of the analytic computations it will obey the same
equations) for (isotropic) quadrics without center, but of a more
general nature and more importantly it will appear as a natural
consequence of the ACPIA.

Once some knowledge about the {\it rigid motion provided by the
Ivory affinity} (RMPIA) is drawn, we shall use this (just as
Bianchi) to provide a simpler proof without using common conjugate
systems.

All computations are local and assumed to be valid on their open
domain of validity without further details; all functions have the
assumed order of differentiability and are assumed to be
invertible, non-zero, etc when required (for all practical
purposes we can assume all functions to be analytic).

\section{Confocal quadrics in canonical form}

Consider the complexified Euclidean space
$$(\mathbb{C}^{n+1},<.,.>),\
<x,y>:=x^Ty,\ |x|^2:=x^Tx,\ x,y\in\mathbb{C}^{n+1}$$ with standard
basis $\{e_j\}_{j=1,...,n+1},\ e_j^Te_k=\del_{jk}$.

Isotropic (null) vectors are those vectors $v$ of length $0\
(|v|^2=0)$; since most vectors are not isotropic we shall call a
vector simply vector and we shall only emphasize isotropic when
the vector is assumed to be isotropic. The same denomination will
apply in other settings: for example we call quadric a
non-degenerate quadric (a quadric projectively equivalent to the
complex unit sphere).

A quadric $x\subset\mathbb{C}^{n+1}$ is given by the equation
$Q(x):=\begin{bmatrix}x\\1\end{bmatrix}^T\begin{bmatrix}A&B\\B^T&C\end{bmatrix}
\begin{bmatrix}x\\1\end{bmatrix}=x^T(Ax+2B)+C=0,\ A=A^T\in\mathbf{M}_{n+1}(\mathbb{C}),\
B\in\mathbb{C}^{n+1},\ C\in\mathbb{C},\
\begin{vmatrix}A&B\\B^T&C\end{vmatrix}\neq 0$.

A metric classification of all (totally real) quadrics in
$\mathbb{C}^{n+1}$ requires the notion of {\it symmetric Jordan}
(SJ) canonical form of a symmetric complex matrix. The symmetric
Jordan blocks are: $J_1:=0=0_{1,1}\in\mathbf{M}_1(\mathbb{C}),\
J_2:=f_1f_1^T\in\mathbf{M}_2(\mathbb{C}),\
J_3:=f_1e_3^T+e_3f_1^T\in\mathbf{M}_3(\mathbb{C}),\ J_4:=f_1\bar
f_2^T+f_2f_2^T+\bar f_2f_1^T\in\mathbf{M}_4(\mathbb{C}),\ J_5:=
f_1\bar f_2^T+f_2e_5^T+e_5f_2^T+\bar
f_2f_1^T\in\mathbf{M}_5(\mathbb{C}),\ J_6:= f_1\bar f_2^T+f_2\bar
f_3^T+f_3f_3^T+\bar f_3f_2^T+\bar
f_2f_1^T\in\mathbf{M}_6(\mathbb{C})$, etc, where
$f_j:=\frac{e_{2j-1}+ie_{2j}}{\sqrt{2}}$ are the standard
isotropic vectors (at least the blocks $J_2,\ J_3$ were known to
the classical geometers). Any symmetric complex matrix can be
brought via conjugation with a complex rotation to the symmetric
Jordan canonical form, that is a matrix block decomposition with
blocks of the form $a_jI_p+J_p$; totally real quadrics are
obtained for eigenvalues $a_j$ of the quadratic part $A$ defining
the quadric being real or coming in complex conjugate pairs $a_j,\
\bar a_j$ with subjacent symmetric Jordan blocks of same dimension
$p$. Just as the usual Jordan block $\sum_{j=1}^pe_je_{j+1}^T$ is
nilpotent with $e_{p+1}$ cyclic vector of order $p$, $J_p$ is
nilpotent with $\bar f_1$ cyclic vector of order $p$, so we can
take square roots of SJ matrices without isotropic kernels
($\sqrt{aI_p+J_p}:=\sqrt{a}\sum_{j=0}^{p-1}(^{\frac{1}{2}}_j)a^{-j}J_p^j,\
a\in\mathbb{C}^*,\ \sqrt{a}:=\sqrt{r}e^{i\theta}$ for
$a=re^{2i\theta},\ 0<r,\ -\pi\le 2\theta<\pi$), two matrices with
same SJ decomposition type (that is $J_p$ is replaced with a
polynomial in $J_p$) commute, etc.

The confocal family $\{x_z\}_{z\in\mathbb{C}}$ of a quadric
$x_0\subset\mathbb{C}^{n+1}$ in canonical form (depending on as
few constants as possible) is given in the projective space
$\mathbb{C}\mathbb{P}^{n+1}$ by the equation
$Q_z(x_z):=\begin{bmatrix}x_z\\1\end{bmatrix}^T(\begin{bmatrix}A&B\\B^T&C\end{bmatrix}^{-1}-z
\begin{bmatrix}I_{n+1}&0\\0^T&0\end{bmatrix})^{-1}\begin{bmatrix}x_z\\1\end{bmatrix}=0$,
where

i) $A=A^T\in\mathbf{GL}_{n+1}(\mathbb{C})$ SJ,
$B=0\in\mathbb{C}^{n+1},\ C=-1$ for {\it quadrics with center}
(QC),

ii) $A=A^T\in\mathbf{M}_{n+1}(\mathbb{C})$ SJ,
$\ker(A)=\mathbb{C}e_{n+1},\ B=-e_{n+1},\ C=0$ for {\it quadrics
without center} (QWC) and

iii) $A=A^T\in\mathbf{M}_{n+1}(\mathbb{C})$ SJ,
$\ker(A)=\mathbb{C}f_1,\ B=-\bar f_1,\ C=0$ for {\it isotropic
quadrics without center} (IQWC).

From the definition one can see that the family of quadrics
confocal to $x_0$ is the adjugate of the pencil generated by the
adjugate of $x_0$ and Cayley's absolute
$C(\infty)\subset\mathbb{C}\mathbb{P}^n$ in the hyperplane at
infinity; since Cayley's absolute encodes the Euclidean structure
of $\mathbb{C}^{n+1}$ (it is the invariant set under rigid motions
and homotheties of
$\mathbb{C}^{n+1}:=\mathbb{CP}^{n+1}\backslash\mathbb{CP}^n$) the
mixed metric-projective character of the confocal family becomes
clear.

For Q(W)C $\mathrm{spec}(A)$ is unambiguous (does not change under
rigid motions) but for IQWC it may change with $(p+1)$-roots of
unity for the block of $f_1$ in $A$ being $J_p$ even under rigid
motions which preserve the canonical form, so it is unambiguous up
to $(p+1)$-roots of unity.

We have the diagonal Q(W)C respectively for
$A=\Sigma_{j=1}^{n+1}a_j^{-1}e_je_j^T,\
A=\Sigma_{j=1}^na_j^{-1}e_je_j^T$; the diagonal IQWC come in
different flavors, according to the block of $f_1:\
A=J_p+\Sigma_{j=p+1}^{n+1}a_j^{-1}e_je_j^T$; in particular if
$A=J_{n+1}$, then $\mathrm{spec}(A)=\{0\}$ is unambiguous. Thus
general quadrics are those for which all eigenvalues have
geometric multiplicity $1$; equivalently each eigenvalue has an
only corresponding SJ block; in this case the quadric also admits
elliptic coordinates.

There are continuous groups of symmetries which preserve the SJ
canonical form for more than one SJ block corresponding to an
eigenvalue, so from a metric point of view a metric classification
according to the elliptic coordinates and continuous symmetries
may be a better one.

With $R_z:=I_{n+1}-zA,\
z\in\mathbb{C}\setminus\mathrm{spec}(A)^{-1}$ the family of
quadrics $\{x_z\}_z$ confocal to $x_0$ is given by
$Q_z(x_z)=x_z^TAR_z^{-1}x_z+2(R_z^{-1}B)^Tx_z+C+zB^TR_z^{-1}B=0$.
For $z\in\mathrm{spec}(A)^{-1}$ we obtain singular confocal
quadrics; those with $z^{-1}$ having geometric multiplicity $1$
admit a singular set which is an $(n-1)$-dimensional quadric
projectively equivalent to $C(\infty)$, so they will play an
important r\^{o}le in the discussion of homographies
$H\in\mathbf{PGL}_{n+1}(\mathbb{C})$ taking a confocal family into
another one, since $H^{-1}(C(\infty)),\ C(\infty)$ respectively
$C(\infty),\ H(C(\infty))$ will suffice to determine each confocal
family.

The Ivory affinity is an affine correspondence between confocal
quadrics and having good metric properties (it may be the reason
why Bianchi calls it {\it affinity} in more than one language): it
is given by $x_z=\sqrt{R_z}x_0+C(z),\
C(z):=-(\frac{1}{2}\int_0^z(\sqrt{R_w})^{-1}dw)B$. Note that
$C(z)=0$ for QC, $=\frac{z}{2}e_{n+1}$ for QWC; for IQWC it is the
Taylor series of $\frac{1}{2}\int_0^z(\sqrt{1-w})^{-1}dw$ at $z=0$
with each monomial $z^{k+1}$ replaced by $z^{k+1}J_p^k\bar f_1$,
where $J_p$ is the block of $f_1$ in $A$ and thus a polynomial of
degree $p$ in $z$. Note
$AC(z)+(I_{n+1}-\sqrt{R_z})B=0=(I_{n+1}+\sqrt{R_z})C(z)+zB$.
Applying $d$ to $Q_z(x_z)=0$ we get $dx_z^TR_z^{-1}(Ax_z+B)=0$, so
the unit normal $N_z$ is proportional to $\hat N_z:=-2\pa_zx_z$.
If $\mathbb{C}^{n+1}\ni x\in x_{z_1},x_{z_2}$, then $\hat
N_{z_j}=R_{z_j}^{-1}(Ax+B)$; using $R_z^{-1}-I_{n+1}=zAR_z^{-1},\
z_1R_{z_1}^{-1}-z_2R_{z_2}^{-1}=(z_1-z_2)R_{z_1}^{-1}R_{z_2}^{-1}$
we get $0=Q_{z_1}(x)-Q_{z_2}(x)=(z_1-z_2)\hat N_{z_1}^T\hat
N_{z_2}$, so two confocal quadrics cut each other orthogonally
(Lam\'{e}). For general quadrics the polynomial equation
$Q_z(x)=0$ has degree $n+1$ in $z$ and it has multiple roots iff
$0=\pa_zQ_z(x)=|\hat N_z|^2$; thus outside the locus of isotropic
normals elliptic coordinates (given by the roots $z_1,...,z_{n+1}$
of the said equation) give a parametrization of $\mathbb{C}^{n+1}$
suited to confocal quadrics.

We have now some classical metric properties of the Ivory
affinity: with $x_0^0,x_0^1\in x_0,\ V_0^1:=x_z^1-x_0^0$, etc the
Ivory Theorem (preservation of length of segments between confocal
quadrics) becomes
$|V_0^1|^2=|x_0^0+x_0^1-C(z)|^2-2(x_0^0)^T(I_{n+1}+\sqrt{R_z})x_0^1+zC=|V_1^0|^2$;
the preservation of lengths of rulings: $w_0^TAw_0=w_0^T\hat
N_0=0,\ w_z=\sqrt{R_z}w_0\Rightarrow
w_z^Tw_z=|w_0|^2-zw_0^TAw_0=|w_0|^2$; the symmetry of the TC:
$(V_0^1)^T\hat N_0^0
=(x_0^0)^TA\sqrt{R_z}x_0^1-B^T(x_z^0+x_z^1-C(z))+C=(V_1^0)^T\hat
N_0^1$; the preservation of angles between segments and rulings:
$(V_0^1)^Tw_0^0+(V_1^0)^Tw_z^0=-z(\hat N_0^0)^Tw_0^0=0$; the
preservation of angles between rulings:
$(w_0^0)^Tw_z^1=(w_0^0)^T\sqrt{R_z}w_0^1=(w_z^0)^Tw_0^1$; the
preservation of angles between polar rulings: $(w_0^0)^TA\hat
w_0^0=0\Rightarrow (w_z^0)^T\hat w_z^0=(w_0^0)^T\hat
w_0^0-z(w_0^0)^TA\hat w_0^0=(w_0^0)^T\hat w_0^0$.

All complex quadrics are affine equivalent to either the unit
sphere $X\subset\mathbb{C}^{n+1},\ |X|^2=1$ or to the equilateral
paraboloid $Z\subset\mathbb{C}^{n+1},\ Z^T(I_{1,n}Z-2e_{n+1})=0$,
so a parametrization with regard to these two quadrics is in
order: $x_0=(\sqrt{A})^{-1}X$ for QC,
$x_0=(\sqrt{A+e_{n+1}e_{n+1}^T})^{-1}Z$ for QWC (for this reason
from a canonical metric point of view (that is we are interested
in a simplest form of $|dx_0|^2$) we should rather require that
$A^{-1}$ or $(A+e_{n+1}e_{n+1}^T)^{-1}$ is SJ).

For IQWC such a parametrization fails because of the isotropic
$\mathrm{ker}(A)$; however the computations between confocal
quadrics involving the Ivory affinity reveal a natural
parametrization of IQWC which is again an affine transformation of
$Z$.

Consider a canonical IQWC $x_0^T(Ax_0-2\bar f_1)=0,\
\ker(A)=\mathbb{C}f_1,\ A=J_p\oplus ...$ SJ. We are looking for a
linear map $L\in\mathbf{GL}_{n+1}(\mathbb{C})$ such that $x_0=LZ$,
equivalently $L^TAL=e^{2a}I_{1,n},\
I_{1,n}:=I_{n+1}-e_{n+1}e_{n+1}^T,\ L^T\bar f_1=e^{2a}e_{n+1}$.
Replacing $L$ with $L(e^{-a}I_{1,n}+e^{-2a}e_{n+1}e_{n+1}^T)$ we
can make $a=0$. Thus $Le_{n+1}=f_1,\ L^T(A+\bar f_1\bar
f_1^T)L=I_{n+1}$, so $L^{-1}=R^T\sqrt{A+\bar f_1\bar f_1^T},\
R^TR=I_{n+1}$ with $Re_{n+1}=\sqrt{A+\bar f_1\bar f_1^T}f_1$ (note
that $Re_{n+1}$ has, as required, length $1$). Once
$R\in\mathbf{O}_{n+1}(\mathbb{C})$ with the above property is
found, $L$ thus defined satisfies $L^T\bar f_1=e_{n+1}$ and thus
$L^TAL=I_{1,n}$. $L$ with the above properties is unique up to
rotations fixing $e_{n+1}$ in its domain and a canonical choice of
$R$ will reveal itself from a SJ canonical form when doing
computations on confocal quadrics. We have $LL^T=(A+\bar f_1\bar
f_1^T)^{-1},\
I_{1,n}L^{-1}\sqrt{R_z}L=I_{1,n}L^{-1}\sqrt{R_z}LI_{1,n}=L^{-1}\sqrt{R_z}L-e_{n+1}\bar
f_1\sqrt{R_z}L=L^TA\sqrt{R_z}L=I_{1,n}\sqrt{I_{n+1}-zL^TA^2L}=:
I_{1,n}\sqrt{R'_z},\ A':=L^TA^2L,\
\ker(A')=\mathbb{C}e_{n+1}\oplus\mathbb{C}L^{-1}(A+\bar f_1\bar
f_1^T)^{-1}f_1 =\mathbb{C}e_{n+1}\oplus\mathbb{C}L^Tf_1$; choose
$R$ which makes $A'$ SJ. Note that we can take for QWC
$L:=(\sqrt{A+e_{n+1}e_{n+1}^T})^{-1},\ A':=A,\
\ker(A')=\mathbb{C}e_{n+1}$, so IQWC can be regarded as metrically
degenerated QWC. Note that
$e_{n+1}^TL^{-1}\sqrt{R_z}L=(-I_{1,n}L^{-1}C(z)+e_{n+1})^T$; this
can be confirmed analytically by differentiating with respect to
$z$ and using $(L^T)^{-1}=AL-Be_{n+1}^T$ and will imply the
symmetry of the TC, but since we have already proved the symmetry
of the TC, we can use this to imply the previous. Thus
$L^{-1}x_z=L^{-1}\sqrt{R_z}LZ+L^{-1}C(z)=I_{1,n}\sqrt{R'_z}Z+
e_{n+1}(-I_{1,n}L^{-1}C(z)+e_{n+1})^TZ+L^{-1}C(z),\
(x_z^1-x_0^0)^T\hat
N_0^0=(L^{-1}x_z^1-Z_0)^T(I_{1,n}Z_0-e_{n+1})=Z_0^TI_{1,n}\sqrt{R'_z}Z_1+
(Z_0+Z_1)^T(I_{1,n}L^{-1}C(z)-e_{n+1})-e_{n+1}^TL^{-1}C(z)$. Note
that for IQWC $|I_{1,n}L^{-1}C(z)|^2=2e_{n+1}^TL^{-1}C(z)$, so
$L^{-1}C(z)$ lies itself on $Z$ (also in this case since $\bar
f_1^TJ_p^k\bar f_1=\del_{k\ p-1}$ we have
$e_{n+1}^TL^{-1}C(z)=\bar
f_1^TC(z)=(^{-\frac{1}{2}}_{p-1})\frac{(-z)^p}{-2p}$, so
$e_{n+1}^T$ picks up the highest power of $z$ in $L^{-1}C(z)$). To
see this we need $0=|L^{-1}C(z)-\bar f_1^TC(z)e_{n+1}|^2-2\bar
f_1^TC(z)=C(z)^T(LL^T)^{-1}C(z)-(\bar f_1^TC(z))^2-2\bar
f_1^TC(z)=C(z)^TAC(z)-2\bar f_1^TC(z)$; using
$AC(z)=(I_{n+1}-\sqrt{R_z})\bar f_1$ and
$(I_{n+1}+\sqrt{R_z})C(z)=z\bar f_1$ it is satisfied.

Note also that as needed later we have
$(L^TL)^{-1}=A'-I_{1,n}L^{-1}Be_{n+1}^T-e_{n+1}(I_{1,n}L^{-1}B)^T+|B|^2e_{n+1}e_{n+1}^T,\
|\hat
N_0|^2=|(L^T)^{-1}(I_{1,n}Z-e_{n+1})|^2=Z^TA'Z+2Z^TI_{1,n}L^{-1}B+|B|^2,\\
(I_{n+1}+\sqrt{R'_z})I_{1,n}L^{-1}C(z)=I_{1,n}L^{-1}(I_{n+1}+\sqrt{R_z})C(z)=-zI_{1,n}L^{-1}B$.

\section{The Bianchi Truths}

We hold these (a-priori not-so-evident) Bianchi Truths to be
self-evident:

{\bf I (existence and inversion of the B\"{a}cklund transformation
for quadrics and the applicability correspondence provided by the
Ivory affinity)}

{\it That any deformation $x^0\subset\mathbb{C}^{2n-1}$ of an
$n$-dimensional sub-manifold $x_0^0\subseteq x_0$
($x_0\subset\mathbb{C}^{n+1}\subset\mathbb{C}^{2n-1}$ being a
quadric) appears as a focal sub-manifold of an
$(\frac{n(n-1)}{2}+1)$-dimensional family of Weingarten
congruences, whose other focal sub-manifolds $x^1=B_{z}(x^0)$ are
applicable, via the Ivory affinity between confocal quadrics, to
sub-manifolds $x_0^1$ in the same quadric $x_0$. The determination
of these sub-manifolds requires the integration of a family of
Riccati equations depending on the parameter $z$ (we ignore for
simplicity the dependence on the initial value data in the
notation $B_z$). Moreover, if we compose the inverse of the rigid
motion provided by the Ivory affinity (RMPIA) $(R_0^1,t_0^1)$ with
the rolling of $x_0^0$ on $x^0$, then we obtain the rolling of
$x_0^1$ on $x^1$ and $x^0$ reveals itself as a $B_z$ transform of
$x^1$;}

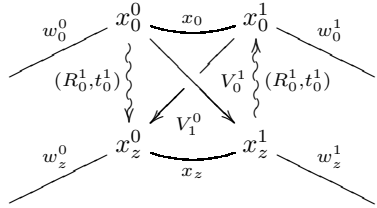
\begin{center}
$\xymatrix@!0{&&x_0^0\ar@{-}[drdr]\ar@/_/@{-}[rr]^{x_0}\ar@{~>}[dd]_{(R_0^1,t_0^1)}&&
x_0^1\ar@{<~}[dd]^{(R_0^1,t_0^1)}&&\\
\ar@{-}[urr]^{w_0^0}&&&\ar[dl]^>>>>{V_1^0}\ar[dr]^>>>>>{V_0^1}&&&\ar@{-}[ull]_{w_0^1}\\&&
x_z^0\ar@{-}'[ur][urur]\ar@/_/@{-}[rr]_{x_z}&&x_z^1&\\\ar@{-}[urr]^{w_z^0}&&&
&&&\ar@{-}[ull]_{w_z^1}}$
\end{center}
\begin{eqnarray}\label{eq:backl}
(R_j,t_j)(x_0^j,dx_0^j):=(R_jx_0^j+t_j,R_jdx_0^j)=(x^j,dx^j),\
(R_j,t_j)\in\mathbf{O}_{2n-1}(\mathbb{C})\ltimes\mathbb{C}^{2n-1},\ j=0,1,\nonumber\\
(R_0^1,t_0^1)=(R_1,t_1)^{-1}(R_0,t_0).
\end{eqnarray}

{\bf II (Bianchi Permutability Theorem)}

{\it That if $x^1=B_{z_1}(x^0),\ x^2=B_{z_2}(x^0)$, then one can
find only by algebraic computations a submanifold
$B_{z_2}(x^1)=x^3=B_{z_1}(x^2)$; thus $B_{z_2}\circ B_{z_1}
=B_{z_1}\circ B_{z_2}$ and once all $B$ transforms of the seed
$x^0$ are found, the $B$ transformation can be iterated using only
algebraic computations;}

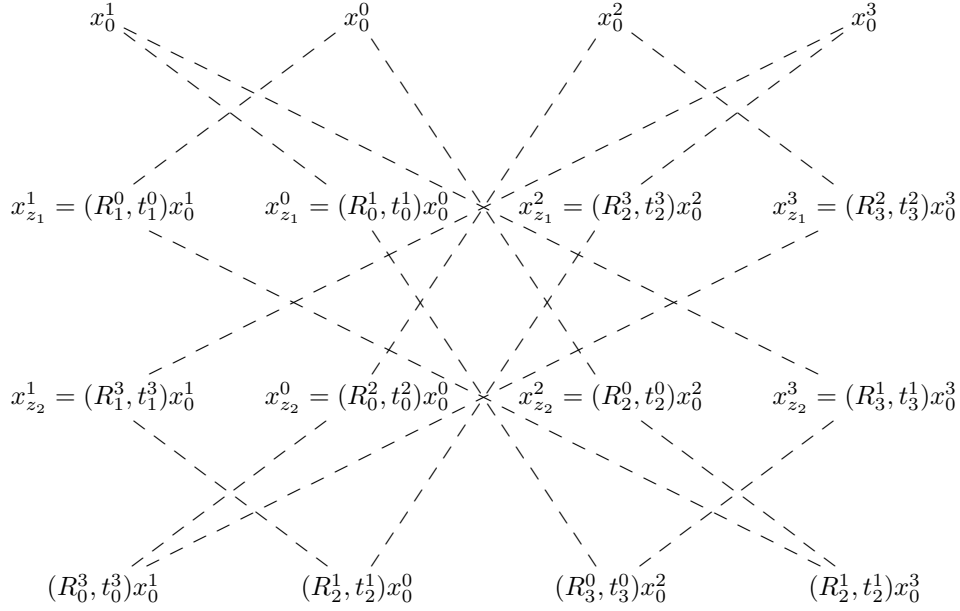
\begin{center}
$\xymatrix@!0{x_0^1\ar@{--}[dddrrrr]\ar@{--}[ddddddrrrrrrrrrrrr]&&&&
x_0^0\ar@{--}[dddllll]\ar@{--}[ddddddrrrr]&&&&
x_0^2\ar@{--}[ddddddllll]\ar@{--}[dddrrrr]&&&&
x_0^3\ar@{--}[ddddddllllllllllll]\ar@{--}[dddllll]\\
\\  \\
x_{z_1}^1=(R_1^0,t_1^0)x_0^1\ar@{--}[ddddddrrrrrrrrrrrr]&&&&
x_{z_1}^0=(R_0^1,t_0^1)x_0^0\ar@{--}[ddddddrrrr]&&&&
x_{z_1}^2=(R_2^3,t_2^3)x_0^2\ar@{--}[ddddddllll]&&&&
x_{z_1}^3=(R_3^2,t_3^2)x_0^3\ar@{--}[ddddddllllllllllll]\\
\\  \\
x_{z_2}^1=(R_1^3,t_1^3)x_0^1\ar@{--}[dddrrrr]&&&&
x_{z_2}^0=(R_0^2,t_0^2)x_0^0\ar@{--}[dddllll]&&&&
x_{z_2}^2=(R_2^0,t_2^0)x_0^2\ar@{--}[dddrrrr]&&&&
x_{z_2}^3=(R_3^1,t_3^1)x_0^3\ar@{--}[dddllll]\\
\\   \\
(R_0^3,t_0^3)x_0^1&&&&(R_2^1,t_2^1)x_0^0&&&&(R_3^0,t_3^0)x_0^2&&&&(R_2^1,t_2^1)x_0^3}$
\end{center}

\begin{eqnarray}\label{eq:cocy}
(R_j,t_j)(x^j,dx^j)=(x_0^j,dx_0^j),\ j=0,...,3,\nonumber\\
(R_j^k,t_j^k)=(R_k,t_k)^{-1}(R_j,t_j),\
(j,k)=(0,1),(0,2),(1,3),(2,3),\\
(R_0^1,t_0^1)(R_2^0,t_2^0)=(R_3^0,t_3^0)=(R_3^1,t_3^1)(R_2^3,t_2^3),\
(R_1^0,t_1^0)(R_3^1,t_3^1)=(R_2^1,t_2^1)=(R_2^0,t_2^0)(R_3^2,t_3^2).\nonumber
\end{eqnarray}

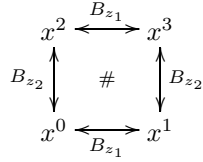
\begin{center}
$\xymatrix{\ar@{}[dr]|{\#}x^2\ar@{<->}[d]_{B_{z_2}}\ar@{<->}[r]^{B_{z_1}}&x^3\ar@{<->}[d]^{B_{z_2}}\\
x^0\ar@{<->}[r]_{B_{z_1}}&x^1}$
\end{center}

{\bf III (existence of $3$-M\"{o}bius moving configurations
$\mathcal{M}_3$)}

{\it That if $x^1=B_{z_1}(x^0),\ x^2=B_{z_2}(x^0),\
x^4=B_{z_3}(x^0)$ and by use of the Bianchi Permutability Theorem
one finds $B_{z_3}(x^2)=x^6=B_{z_2}(x^4),\
B_{z_1}(x^4)=x^5=B_{z_3}(x^1),\ B_{z_2}(x^1)=x^3=B_{z_1}(x^2),\
B_{z_3}(x^3)=x'^7=B_{z_2}(x^5),\ B_{z_1}(x^6)=x''^7=B_{z_3}(x^3),\
B_{z_2}(x^5)=x'''^7=B_{z_1}(x^6)$, then $x'^7=x''^7=x'''^7=:x^7$;
thus once all $B$ transforms of the seed $x^0$ are found, the $B$
transformation can be further iterated using only algebraic
computations;}

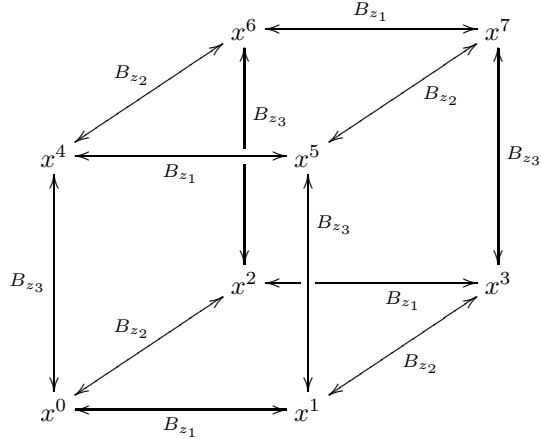
\begin{center}
$\xymatrix@!0{&&&x^6\ar@{<->}[rrrr]^{B_{z_1}}&&&&x^7\\
&&&&&&&\\
x^4\ar@{<->}[uurrr]^{B_{z_2}} \ar@{<->}[rrrr]_{B_{z_1}}&&&&
x^5\ar@{<->}[uurrr]_>>>>>>>>>{B_{z_2}}&&&\\
&&&&&&&\\
&&&x^2\ar@{<->}'[r][rrrr]_{B_{z_1}}
\ar@{<->}'[uu][uuuu]_<<<<<{B_{z_3}}&&&&
x^3\ar@{<->}[uuuu]_{B_{z_3}}\\
&&&&&&&\\
x^0\ar@{<->}[rrrr]_{B_{z_1}}\ar@{<->}[uurrr]^{B_{z_2}}
\ar@{<->}[uuuu]^{B_{z_3}}&&&& x^1\ar@{<->}[uurrr]_{B_{z_2}}
\ar@{<->}[uuuu]_<<<<<<<<<<<<<<<<<<<<<<{B_{z_3}}&&&}$
\end{center}

{\bf IV (Hazzidakis transformation)}

{\it That if an $n$-dimensional sub-manifold
$x^0\subset\mathbb{C}^{2n-1}$ is applicable to a sub-manifold
$x_0^0\subseteq x_0$ and the homography
$H\in\mathbf{PGL}_{n+1}(\mathbb{C})$ takes the confocal family
$x_z$ to another confocal family $\ti x_{\ti z},\ \ti z=\ti z(z),\
\ti z(0)=0$, then one infinitesimally knows a sub-manifold $\ti
x^0=H(x^0)$ (that is one knows the first and second fundamental
forms), called the Hazzidakis (H) transform of $x^0$ and
applicable to a sub-manifold $\ti x_0^0\subseteq\ti x_0$. Moreover
the H transformation commutes with the B transformation ($H\circ
B_z=B_{\ti z}\circ H$) and the $B_{\ti z}(\ti x^0)$ transforms can
be algebraically recovered from the knowledge of $\ti x^0$ and
$B_z(x^0)$}.

Keeping an eye on

{\bf 0 (The Method of Archimedes)}

{\it '... certain things first became clear to me by a mechanical
method, although they had to be proved by geometry afterwards
because their investigation by the said method did not furnish an
actual proof. But it is of course easier, when we have previously
acquired, by the method, some knowledge of the questions, to
supply the proof than it is to find it without any previous
knowledge'}

we ask ourselves what happens when the {\it 'mechanical method'}
in question is the rolling. In this case {\it 'first'} means that
the seed $x^0$ is in the particular position of actually
coinciding with $x_0^0\subseteq
x_0\subset\mathbb{C}^{n+1}\subset\mathbb{C}^{2n-1}$: for $n=2$ the
leaves become rulings on a confocal quadric (Bianchi and Lie; in
fact it is this small observation of Lie's for the pseudo-sphere
that proved to be the essential tool used by Bianchi to develop
the theory of deformations of quadrics, but the classical
geometers did not investigate what becomes of Lie's powerful {\it
facets} (pairs of points and planes passing through those points)
approach when the leaf degenerates from a surface to a ruling).

Thus the natural conjecture appears that the B transforms $x^1$
(leaves) of $x^0$ (seed) are rulings on quadrics $x_z$ confocal to
$x_0$ when $x^0=x_0^0\subset x_0$. An a-priori argument in favor
of rulings on quadrics $x_z$ confocal to $x_0$ being the leaves
when $x^0=x_0^0$ is the fact that in this case the leaves should
be sub-manifolds of $x_z$; thus their tangent spaces should be
contained both in the tangent bundle of $x_z$ and in the
distribution of facets; since facets intersect tangent spaces of
$x_z$ only along rulings of $x_z$ the leaves should be rulings.
Also rulings preserve their length under the Ivory affinity, so
they should be naturally chosen by the RMPIA (the discretization
of the infinitesimal version of the ACPIA).

Due to the symmetries (in the normal bundle) of the rolling we get
the
$\mathbf{O}_{n-1}(\mathbb{C})\times\mathbf{O}_{n-1}(\mathbb{C})$
symmetries of the facets in the tangency configuration and due to
the initial value data of the differential system subjacent to the
B transformation we get the fact that for the deformation problem
of $n$-dimensional quadrics in $\mathbb{C}^{2n-1}$ no functional
information is allowed in the normal bundle.

Since by rolling the position of the leaf $x^1$ relative to the
seed $x^0$ is the same as the position of $x_z^1$ relative to
$x_0^0$, we get the fact that the tangential part of $dx^1$
(relative to $x^0$) is the same as the tangential part of $dx_z^1$
(relative to $x_0^0$). Thus with $N^0:=[N_{n+1}^0\ \ ...\ \
N_{2n-1}^0]$ orthonormal normal frame of $x^0$ and $N_0^0$ unit
normal of $x_0^0$ we have
$|dx^1|^2=|(dx^1)^{\top}|^2+|(dN^0)^T(x^1-x^0)|^2,\
|dx_z^1|^2=|(dx_z^1)^{\top}|^2+|(dN_0^0)^T(x_z^1-x_0^0)|^2$ and
the ACPIA $|dx^1|^2=|dx_0^1|^2$ becomes
\begin{eqnarray}\label{eq:fund}
|dx_0^1|^2-|dx_z^1|^2=|\begin{bmatrix}-i(dN_0^0)^T(x_z^1-x_0^0)\\
-(dN^0)^T(x^1-x^0)\end{bmatrix}|^2
\end{eqnarray} (we augment the
$(n-1)$-column $1$-form $-(dN^0)^T(x^1-x^0)$ with the entry
$-i(dN_0^0)^T(x_z^1-x_0^0)$ on the first position to get an
$n$-column $1$-form; thus the second fundamental forms are
naturally joined).

Equation (\ref{eq:fund}) will turn out to be essential, since by
the metric properties of the Ivory affinity the {\it left hand
side} (lhs) will be an $n$-dimensional symmetric quadratic form
$|\om_1|^2$ which is a square for proper choice of coordinates,
just like the {\it right hand side} (rhs), so one obtains an
equality $\om_1=R\om_0$ of $n$-column $1$-forms (involving an
a-priori arbitrary rotation matrix
$R\subset\mathbf{O}_n(\mathbb{C})$).

However, at least for the first draft a coordinate-independent
approach seems to stop here: to continue two sets of coordinates
must be used in order to take full advantage of (\ref{eq:fund}): a
set of coordinates to express in a simplest form the lhs and other
for the rhs of (\ref{eq:fund}), while keeping in mind meaningful
formulae for changes of coordinates (the change of Christoffel
symbols will turn out to be the main problem). For the lhs we need
to take into consideration only two cases of quadrics: the
equilateral paraboloid (when the parametrization is the one that
realizes $x_0$ as a graph) and the unit sphere (when the
parametrization is given by the stereo-graphical projection; it is
the projective transformation of the graph coordinates on the
equilateral paraboloid under the projective transformation which
takes the point at infinity on the $e_{n+1}$-axis to the north
pole); these parametrizations have the advantage of providing an
isothermic-conjugate system on $x_0$. For the rhs we shall use
conjugate systems common to both $x_0^0$ and $x^0$, since they
will considerably simplify the computations and will naturally fit
for the deformation problem for higher dimensional quadrics, just
as for $2$-dimensional quadrics.

The investigation of the first three Bianchi Truths respectively
boil down when $x^0=x_0^0$ to the investigation of the {\it first,
second and third iterations of the tangency configuration} (TC,
SITC, TITC). The Bianchi Truths I - IV become respectively the
Bianchi I - IV Theorems on confocal quadrics (some of the metric
properties of the Ivory affinity between confocal quadrics and
beyond the Ivory Theorem on the preservation under the Ivory
affinity of lengths of segments between confocal quadrics were
already known to other authors; for example Henrici's construction
of the articulated hyperbolic paraboloid uses preservation of
lengths of rulings under the Ivory affinity).

Note that for $n=2$ the Bianchi Truth III uses a theorem of
Menelaus (in itself a co-cycle theorem), which is equivalent to
the infinitesimal associativity of a loop group action.

Since the dimensionality of the space of leaves should be
independent of the shape of the seed (and thus should equal the
dimensionality of the space of rulings), the two coincide
($2n-3=\frac{n(n-1)}{2}$) without requiring multiplicities of
rulings only for $n=2,3$.

Note that the SITC for an edge of the Bianchi quadrilateral being
infinitesimal infinitesimally describes the B transformation; thus
the SITC encodes all necessary and sufficient information needed
to prove the Bianchi Truth I. A similar statement holds for the
BPT and the TITC, so the first three iterations of the tangency
configuration contain all necessary algebraic information needed
to develop the theory of deformations of quadrics, as expected (by
discretization each iteration corresponds to a derivative, so each
iteration encodes respectively the information in the tangent
space, the {\it Gau\ss-Weingarten} (GW) equations and the {\it
Gau\ss-Codazzi-Mainardi(-Peterson)-Ricci} (G-CMP-R) equations).

While a simple conjugation trick in the BPT provides a tool to
obtain totally real deformations of the same metric type and in
the same Lorentz surrounding space
$\mathbb{R}^{m}\times(i\mathbb{R})^{2n-1-m},\ 0\le m\le 2n-1$ as
the seed from totally real seed without worrying about the
intermediary leaves, a full discussion of the cases when both the
leaf and the seed are totally real (thus they must be situated in
the same Lorentz space) is not completed even for $n=2$.

\section{The deformation problem for quadrics via common conjugate systems}

First we shall recall the notion of deformations of quadrics with
common conjugate system and non-degenerate joined second
fundamental forms, appearing in one of our previous notes
concerning Peterson's deformations of higher dimensional quadrics.

Consider the complexified Euclidean space
$$(\mathbb{C}^m,<.,.>),\
<x,y>:=x^Ty,\ |x|^2:=x^Tx,\ x,y\in\mathbb{C}^m$$ with standard
basis $\{e_j\}_{j=1,...,m},\ e_j^Te_k=\del_{jk}$ (we shall use
$m=n,n+1,2n-1$).

We shall always have Latin indices $j,k,l,...\in\{1,...,n\}$
(including for differentiating respectively with respect to $u^j,\
u^k,\ u^l,...$), Greek ones $\al,\be,\ga,...\in\{n+1,...,2n-1\}$
for $x\subset\mathbb{C}^{2n-1}$ (or $\al=0$ for
$x_0\subset\mathbb{C}^{n+1}$) and mute summation for upper and
lower indices when clear from the context; also we shall preserve
the classical notation $d^2$ for the tensorial (symmetric) second
derivative and we shall use $d\wedge$ for the exterior
(antisymmetric) derivative; thus $d\wedge d=0$.

For $n\ge 3$ consider the $n$-dimensional sub-manifold
$$x=x(u^1,u^2,...,u^n)\subset\mathbb{C}^{n+p},\ du^1\wedge
du^2\wedge...\wedge du^n\neq 0,\ p=1,n-1$$ such that the tangent
space at any point of $x$ is not isotropic (the scalar product
induced on it by the Euclidean one on $\mathbb{C}^{n+p}$ is not
degenerate; this assures the existence of orthonormal normal
frames). We have the normal frame $N:=[N_{n+1}\ \ N_{n+p}],\
N^TN=I_p$, the first $|dx|^2=g_{jk}du^j\odot du^k$ and second
$d^2x^TN=[h_{jk}^{n+1}du^j\odot du^k\ \ ...\ \
h_{jk}^{n+p}du^j\odot du^k]$ fundamental forms, the Christoffel
symbols
$\Ga_{jk}^l=\frac{g^{lm}}{2}[(g_{jm})_k+(g_{km})_j-(g_{jk})_m]$,
the Riemann curvature
$R_{jmkl}=g_{mp}R^p_{jkl}=g_{mp}[(\Ga_{jk}^p)_l-(\Ga_{jl}^p)_k+\Ga_{jk}^q\Ga_{ql}^p
-\Ga_{jl}^q\Ga_{qk}^p]$ tensor, the normal connection
$N^TdN=\{n^{\al}_{\be j}du^j\}_{\al,\be=n+1,...,n+p},\
n^{\al}_{\be j}=-n^{\be}_{\al j}$ and the curvature $r_{\al
jk}^{\be}=(n_{\al j}^{\be})_k-(n_{\al k}^{\be})_j+n_{\al
j}^{\ga}n_{\ga k}^{\be}-n_{\al k}^{\ga}n_{\ga j}^{\be}$ tensor of
the normal bundle.

We have the {\it Gau\ss-Weingarten} (GW) equations

$$x_{jk}=\Ga_{jk}^lx_l+h_{jk}^{\al}N_{\al},\ (N_{\al})_j=
-h_{jk}^{\al}g^{kl}x_l+n^{\be}_{j\al}N_{\be}$$ and their
integrability conditions $x_{jkl}=x_{jlk},\
(N_{\al})_{jk}=(N_{\al})_{kj}$, from where one obtains by taking
the tangential and normal components (using
$-(g^{jk})_l=g^{jm}\Ga_{ml}^k+g^{km}\Ga_{ml}^j$ and the GW
equations themselves) the {\it
Gau\ss-Codazzi-Mainardi(-Peterson)-Ricci} (G-CMP-R) equations

$$R_{jmkl}=\sum_{\al}(h_{jk}^{\al}h_{lm}^{\al}-h_{jl}^{\al}h_{km}^{\al}),\
(h_{jk}^{\al})_l-(h_{jl}^{\al})_k+\Ga_{jk}^mh_{ml}^{\al}-\Ga_{jl}^mh_{mk}^{\al}
+h_{jk}^{\be}n_{\be l}^{\al}-h_{jl}^{\be}n_{\be k}^{\al}=0,$$
$$r_{\al jk}^{\be}=h_{jl}^{\al}g^{lm}h_{mk}^{\be}-h_{kl}^{\al}g^{lm}h_{mj}^{\be}.$$
If we have conjugate system $h_{jk}^{\al}=:\del_{jk}h_j^{\al}$,
then the above equations become:
\begin{eqnarray}\label{eq:riem}
R_{jkjk}=-R_{jkkj}=\sum_{\al}h_j^{\al}h_k^{\al},\
(h_j^{\al})_k=\Ga_{jk}^jh_j^{\al}-\Ga_{jj}^kh_k^{\al}-h_j^{\be}\eta_{\be
k}^{\al},\ j\neq k,\
R_{jklm}=0\ \mathrm{otherwise},\nonumber\\
\Ga_{jk}^lh_l^{\al}=\Ga_{jl}^kh_k^{\al},\ j,k,l\
\mathrm{distinct},\ r_{\al
jk}^{\be}=(h_j^{\al}h_k^{\be}-h_j^{\be}h_k^{\al})g^{jk}.
\end{eqnarray}
In particular for lines of curvature parametrization
($g_{jk}=\del_{jk}g_{jk}$) we have flat normal bundle, so one can
choose up to multiplication on the right by a constant matrix
$\in\mathbf{O}_p(\mathbb{C})$ normal frame $N$ with zero normal
connection $N^TdN=0$.

Consider now $x\subset\mathbb{C}^{2n-1}$ deformation of
$x_0\subset\mathbb{C}^{n+1}$ (that is $|dx_0|^2=|dx|^2$) with
common conjugate system $(u^1,...,u^n)$ and non-degenerate joined
second fundamental forms (that is $[d^2x_0^TN_0\ \ d^2x^TN]$ is a
symmetric quadratic $\mathbb{C}^n$-valued form which contains only
$(du^j)^2$ terms for $N_0$ normal field of $x_0$ and $N=[N_{n+1}\
\ ...\ \ N_{2n-1}]$ normal frame of $x$ and the dimension $n$
cannot be lowered for $(u^1,...,u^n)$ in an open dense set).

By the argument of Cartan's reduction of exteriorly orthogonal
forms to the canonical form such coordinates $\{u^j\}_j$ exist for
real deformations $\subset\mathbb{R}^{2n-1}$ of imaginary quadrics
$\subset\mathbb{R}^n\times i\mathbb{R}$ (for such cases also all
computations in the deformation problem will be real; see
\cite{B}). However, since we shall derive completely integrable
differential systems (systems in involution) for the deformation
problem for quadrics, the dimensionality of the space of
deformations of quadrics remains the same (namely solution
depending on $n(n-1)$ functions of one variable) in the complex
setting also (the Cartan characters remain the same).

The common conjugate system with non-degenerate joined second
fundamental forms assumption amounts to the vectors $h_j:=[ih_j^0\
\ h_j^{n+1}\ \ ...\ \ h_j^{2n-1}]^T$ being linearly independent.
From the Gau\ss\ equations we obtain
$h_j^0h_k^0=R_{jkjk}=\sum_{\al}h_j^{\al}h_k^{\al},\ j\neq
k\Leftrightarrow h_j^Th_k=\del_{jk}|h_j|^2$; thus the vectors
$h_j\subset\mathbb{C}^n$ are further orthogonal, which prevents
them from being isotropic (should one of them be isotropic, by a
rotation of $\mathbb{C}^n$ one can make it $f_1$
($\mathbf{O}_n(\mathbb{C})$ acts transitively on the isotropic
cone) and after subtracting suitable multiples of $f_1$ from the
remaining ones by another rotation of $\mathbb{C}^n$ the remaining
ones linear combinations of $e_3,...,e_n$, so we would have $n-1$
linearly independent orthogonal vectors in $\mathbb{C}^{n-2}$, a
contradiction), so $\mathbf{a}_j:=|h_j|\neq 0,\
h_j=:\mathbf{a}_jv_j,\ [v_1\ \ ...\ \
v_n]\subset\mathbf{O}_n(\mathbb{C})$.

Note that the linear element must satisfy the condition
$$\Ga_{jk}^l=0,\ j,k,l\ \mathrm{distinct}$$ and deformations
$x\subset\mathbb{C}^{2n-1}$ of $x_0\subset\mathbb{C}^{n+1}$ with
common conjugate system and non-degenerate joined second
fundamental forms are in bijective correspondence with solutions
of the differential system
\begin{eqnarray}\label{eq:syst0}
(\log\mathbf{a}_j)_k=\Ga_{jk}^j,\ j\neq k,\
\sum_j\frac{(h_j^0)^2}{\mathbf{a}_j^2}+1=0.
\end{eqnarray}
Once a solution of this system is known, one finds the second
fundamental form of $x$ and then one finds $x$ by the integration
of a Ricatti equation and quadratures (the
Gau\ss-Bonnet(-Peterson) Theorem).

Note that with $\ga_{jk}:=\Ga_{jj}^k\frac{h_k^0}{h_j^0},\ j\neq k$
the condition
\begin{eqnarray}\label{eq:comint}
(\ga_{jk})_j=(\ga_{kj})_k=-2\ga_{jk}\ga_{kj},\
(\ga_{lj})_k=2(\ga_{lj}\ga_{lk}-\ga_{lk}\ga_{kj}-\ga_{lj}\ga_{jk}),\
j,\ k,\ l\ \mathrm{distinct}
\end{eqnarray}
is necessary and sufficient to get a maximal ($(n-1)$-dimensional)
family of Peterson's deformations $x$ of $x_0$ with common
conjugate system $(u^1,...,u^n)$ and non-degenerate joined second
fundamental forms; it is invariant under transformations of the
variables $u^j$ into themselves, so it is a statement about curves
of coordinates $(u^1,...,\hat u^j,...,u^n)=$ct, $j=1,...,n$.

\subsection{Deformations of (isotropic) quadrics without center}\noindent

\noindent For the specific computations of deformations of
quadrics we shall use the convention
$\mathbb{C}^n\subset\mathbb{C}^{n+1}$ with $0$ on the
$(n+1)^{\mathrm{th}}$ component; thus for example we can multiply
$(n+1,n+1)$-matrices with $n$-column vectors and similarly one can
extend $(n,n)$ matrices to $(n+1,n+1)$ matrices with zeroes on the
last column and row. The converse is also valid: an $(n+1,n+1)$
matrix with zeroes on the last column and row (or multiplied on
the left with an $n$-row vector and on the right with an
$n$-column vector) will be considered as an $(n,n)$-matrix.

With $V:=\sum_{k=1}^nv^ke_k=[v^1\ \ ...\ \ v^n]^T$ consider the
complex equilateral paraboloid
$Z=Z(v^1,...,v^n)=V+\frac{|V|^2}{2}e_{n+1}$.

We have the (I)QWC $x_0:=LZ,\ L\in\mathbf{GL}_{n+1}(\mathbb{C})$
(recall $L:=(\sqrt{A+e_{n+1}e_{n+1}^T})^{-1},\
\ker(A)=\mathbb{C}e_{n+1}, A$ SJ, $B=-e_{n+1}$ for QWC and
$Le_{n+1}=f_1,\ L^T(A+\bar f_1\bar f_1^T)L=I_{n+1},\ A':=L^TA^2L$
SJ for $\ker(A)=\mathbb{C}f_1, A$ SJ, $B=-\bar f_1$ in the case of
IQWC) with linear element, unit normal, second fundamental form
and Christoffel symbols
$|dx_0|^2=dV^TL^TLdV+(V^TdV)^2|Le_{n+1}|^2+2(Le_{n+1})^TLdV(V^TdV),\
N_0=\frac{(L^T)^{-1}V+B}{\sqrt{H}},\
N_0^Td^2x_0=-\frac{|dV|^2}{\sqrt{H}},\
H:=|(L^T)^{-1}V-B|^2=V^TA'V+2V^TL^{-1}B+|B|^2,\ \ti\Ga_{jk}^l=0,\
j\neq k,\ \ti\Ga_{jj}^k=\frac{\pa\log\sqrt{H}}{\pa v^k}$. Note
that we have a distinguished tangent vector field
$\mathcal{V}_0:=\sum_{k=1}^n\frac{\pa\log\sqrt{H}}{\pa
v^k}x_{0v^k}=(x_{0v^jv^j})^{\top}=-\frac{1}{\sqrt{H}}N_0+Le_{n+1}$;
it has the property
$dx_{0v^j}=\mathcal{V}_0dv^j-N_0dN_0^Tx_{0v^j}$.

The condition (\ref{eq:comint}) that $(v^1,...,v^n)$ are common to
a Peterson's $(n-1)$-dimensional family of deformations becomes
$\frac{\pa^2\log\sqrt{H}}{\pa v^j\pa
v^k}=-2\frac{\pa\log\sqrt{H}}{\pa v^j}\frac{\pa\log\sqrt{H}}{\pa
v^k},\ j\neq k\Leftrightarrow e_j^TA'e_k=0,\ j\neq k$, so the
first $(n,n)$ entries of $A'$ must be a diagonal matrix.

Because $(v^1,...,v^n)$ are isothermal-conjugate and
$(u^1,..,u^n)$ are conjugate on $x_0$, the Jacobian
$\frac{\pa(v^1,...,v^n)}{\pa(u^1,...,u^n)}$ has orthogonal
columns, so with $\la_j:=|\frac{\pa V}{\pa u^j}|,\ \La:=[\la_1\ \
...\ \ \la_n]^T,\ \del:=\mathrm{diag}[du^1\ \ ...\ \ du^n]$ we
have $dV=R\del\La,\ R\subset\mathbf{O}_n(\mathbb{C})$. Multiplying
the formula for the change of Christoffel symbols $\frac{\pa
v^c}{\pa u^l}\Ga_{jk}^l=\frac{\pa^2v^c}{\pa u^j\pa u^k}+\frac{\pa
v^a}{\pa u^j}\frac{\pa v^b}{\pa u^k}\ti
\Ga_{ab}^c=\frac{\pa^2v^c}{\pa u^j\pa
u^k}+\la_j^2\del_{jk}\frac{\pa\log\sqrt{H}}{\pa v^c}$ on the left
with $\frac{\pa v^c}{\pa u^p}$ and summing after $c$ we obtain
$\Ga_{jk}^p=\la_p^{-2}(\sum_c\frac{\pa v^c}{\pa
u^p}\frac{\pa^2v^c}{\pa u^j\pa
u^k}+\la_j^2\del_{jk}(\log\sqrt{H})_p)=\del_{pk}(\log\la_k)_j
+\del_{jk}\frac{\la_j^2}{\la_p^2}(\log\frac{\sqrt{H}}{\la_j})_p
+\del_{pj}(\log\la_j)_k$, so $\Ga_{jk}^j=(\log\la_j)_k,\
\Ga_{jj}^k=\frac{\la_j^2}{\la_k^2}(\log\frac{\sqrt{H}}{\la_j})_k,\
j\neq k,\ \Ga_{jj}^j=(\log(\la_j\sqrt{H}))_j$. We have
$h_j^0=-\frac{\la_j^2}{\sqrt{H}}$; since
$(\log\la_j)_k=\Ga_{jk}^j=(\log\mathbf{a}_j)_k,\ j\neq k$ we get
$\la_j=\phi_j(u^j)\mathbf{a}_j$; after a change of the $u^j$
variable into itself we can make $\la_j=\mathbf{a}_j,\ j=1,...,n$,
so from (\ref{eq:syst0}) $|\La|^2=-H$ and
$\La^Td\La=-dV^T(A'V+L^{-1}B)=-\La^T\del R^T(A'V+L^{-1}B)$.

Imposing the compatibility condition $R^Td\wedge$ on $dV=R\del\La$
we get $R^TdR\wedge\del\La-\del\wedge d\La=0$, or
$(\la_j)_k=e_j^TR^TR_je_k\la_k,\ j\neq k,\
e_l^T\frac{R^TR_j}{\la_j}e_k=e_l^T\frac{R^TR_k}{\la_k}e_j,\ j,k,l$
distinct. Now by the standard Cartan trick
$-e_k^T\frac{R^TR_l}{\la_l}e_j=-e_k^T\frac{R^TR_j}{\la_j}e_l=e_l^T\frac{R^TR_j}{\la_j}e_k
=e_l^T\frac{R^TR_k}{\la_k}e_j=-e_j^T\frac{R^TR_k}{\la_k}e_l=-e_j^T\frac{R^TR_l}{\la_l}e_k
=e_k^T\frac{R^TR_l}{\la_l}e_j$ for $j,k,l$ distinct, so
$e_j^TR^TR_ke_l=0$ for $j,k,l$ distinct. Keeping account of the
prime integral property $\La^T[d\La+\del R^T(A'V+L^{-1}B)]=0$ and
with $\om:=\sum_{j=1}^n(e_je_j^TR^TR_j\del+\del
R^TR_je_je_j^T)=-\om^T$ we have $d\La=\om\La-\del
R^T(A'V+L^{-1}B)$ and $d\wedge d\La=0$ becomes the differential
system
\begin{eqnarray}\label{eq:defqwc}
d\wedge\om-\om\wedge\om=-\del R^TA'R\wedge\del,\
\om\wedge\del-\del\wedge R^TdR=0
\Leftrightarrow\nonumber\\
e_j^T[(R^TR_j)_j-(R^TR_k)_k-\sum_lR^TR_le_le_l^TR^TR_l
+R^TA'R]e_k=0,\ j\neq k,\nonumber\\
e_j^TR^TR_le_k=0\ \mathrm{for}\ j,k,l\ \mathrm{distinct},\
R\subset\mathbf{O}_n(\mathbb{C})
\end{eqnarray}
in involution (that is no further conditions appear if one imposes
$d\wedge$ conditions and one uses the equations of the system
itself) as the compatibility condition for the completely
integrable linear system
\begin{eqnarray}\label{eq:systqwc}
dV=R\del\La,\
d\La=\om\La-\del R^T(A'V+L^{-1}B),\ \La^T\La=-(V^TA'V+2V^TL^{-1}B+|B|^2),\nonumber\\
\om:=\sum_{j=1}^n(e_je_j^TR^TR_j\del+\del R^TR_je_je_j^T).
\end{eqnarray}
The differential system (\ref{eq:defqwc}) is similar in some
aspects to that of Terng's GSGE; in fact for $n=2,\
A:=\mathrm{diag}[a_1^{-1}\ \ a_2^{-1}\ \ 0]$ and a normalization
$a_1^{-1}-a_2^{-1}=1$ it is the sine-Gordon equation, but the
quadratic dependence on $R$ in the first equation for Terng's GSGE
is only along the first column $Re_1$ since a similar phenomenon
appears in the definition of $\om$.

Note that given a solution of (\ref{eq:defqwc}) (\ref{eq:systqwc})
will produce an $(n-1)$-dimensional family of deformations
$x\subset\mathbb{C}^{2n-1}$ (the prime integral property removes a
dimension and translations in $u^j$ another $n$ from the original
$2n$-dimensional space of solutions).

One can probably directly analytically develop the B
transformation and its BPT of (\ref{eq:defqwc}), but we shall not
insist on it now, since it will appear naturally at the level of
the geometric picture from the B transformation of quadrics. Note
however that one can deduce at this point that for the first
$(n,n)$ entries of $A'$ being a diagonal matrix the $0$-soliton
$R=I_n$ of (\ref{eq:defqwc}) will produce Peterson's deformations
of quadrics and $\La$ with some entries zero will produce
degenerate solutions (including $0$-solitons of paraboloids), so
these Peterson's deformations of quadrics will be amenable to
explicit computations of their B transforms (at each iteration of
the B transformation the $(n-1)$-dimensional family of solutions
obtained form the passage of solutions of (\ref{eq:defqwc}) to
solutions of (\ref{eq:systqwc}) will be heirs of the original
$(n-1)$-dimensional family of Peterson's deformations of
quadrics), but this is not our interest right now.

\subsection{Deformations of quadrics with center}\noindent

\noindent

With $V:=\sum_{k=1}^nv^ke_k=[v^1\ \ ...\ \ v^n]^T$ consider the
complex unit sphere
$X=X(v^1,...,v^n)=\frac{2V+(|V|^2-1)e_{n+1}}{|V|^2+1}$. We have
$dX=2\frac{dV+V^TdV(e_{n+1}-X)}{|V|^2+1}$, so
$|dX|^2=\frac{4|dV|^2}{(|V|^2+1)^2}$.

With $A=A^T\in\mathbf{GL}_{n+1}(\mathbb{C})$ SJ, $H:=X^TAX$ we
have the QC $x_0=(\sqrt{A})^{-1}X$ with linear element, unit
normal, second fundamental form and Christoffel symbols
$|dx_0|^2=dX^TA^{-1}dX,\ N_0=\frac{\sqrt{A}X}{\sqrt{H}},\
N_0^Td^2x_0=-\frac{|dX|^2}{\sqrt{H}},\ \ti\Ga_{jk}^l=0,\ j,k,l$
distinct, $\ti\Ga_{jk}^j=-\frac{\pa\log(|V|^2+1)}{\pa v^k},\
\ti\Ga_{jj}^k= \frac{\pa\log(\sqrt{H}(|V|^2+1))}{\pa v^k},\ j\neq
k,\ \ti\Ga_{jj}^j=\frac{\pa\log\frac{\sqrt{H}}{|V|^2+1}}{\pa
v^j}$. Note that we have a distinguished tangent vector field
$\mathcal{V}_0:=\sum_{k=1}^n\frac{\pa\log\sqrt{H}}{\pa
v^k}x_{0v^k}\\=\frac{4}{(|V|^2+1)^2}(\frac{N_0}{\sqrt{H}}-x_0)$;
it has the property $dx_{0v^j}=-\frac{\pa\log(|V|^2+1)}{\pa
v^j}dx_0-d\log(|V|^2+1)x_{0v^j}+(\mathcal{V}_0+\sum_{k=1}^n\frac{\pa\log(|V|^2+1)}{\pa
v^k}x_{0v^k})dv^j-N_0dN_0^Tx_{0v^j}$.

Because $(v^1,...,v^n)$ are isothermic-conjugate and
$(u^1,..,u^n)$ are conjugate on $x_0$, the Jacobian
$\frac{\pa(v^1,...,v^n)}{\pa(u^1,...,u^n)}$ has orthogonal
columns, so with $\la_j:=|\frac{\pa V}{\pa u^j}|,\ \La:=[\la_1\ \
...\ \ \la_n]^T,\ \del:=\mathrm{diag}[du^1\ \ ...\ \ du^n]$ we
have $dV=R\del\La,\ R\subset\mathbf{O}_n(\mathbb{C})$. Multiplying
the formula for the change of Christoffel symbols $\frac{\pa
v^c}{\pa u^l}\Ga_{jk}^l=\frac{\pa^2v^c}{\pa u^j\pa u^k}+\frac{\pa
v^a}{\pa u^j}\frac{\pa v^b}{\pa u^k}\ti \Ga_{ab}^c$ on the left
with $\frac{\pa v^c}{\pa u^p}$ and summing after $c$ we obtain
$\Ga_{jk}^p=\frac{1}{\la_p^2}(\sum_c\frac{\pa v^c}{\pa
u^p}\frac{\pa^2v^c}{\pa u^j\pa u^k}+\sum_c\frac{\pa v^c}{\pa
u^p}\frac{\pa v^a}{\pa u^j}\frac{\pa v^b}{\pa u^k}\ti
\Ga_{ab}^c)=\del_{pk}(\log\frac{\la_k}{|V|^2+1})_j+\del_{pj}(\log\frac{\la_j}{|V|^2+1})_k
+\del_{jk}\frac{\la_j^2}{\la_p^2}(\log\frac{\sqrt{H}(|V|^2+1)}{\la_j})_p$,
so $\Ga_{jk}^j=(\log\frac{\la_j}{|V|^2+1})_k,\
\Ga_{jj}^k=\frac{\la_j^2}{\la_k^2}(\log\frac{\sqrt{H}(|V|^2+1)}{\la_j})_k,\
j\neq k,\ \Ga_{jj}^j=(\log\frac{\la_j\sqrt{H}}{|V|^2+1})_j$. We
have $h_j^0=\frac{-4\la_j^2}{\sqrt{H}(|V|^2+1)^2}$; since
$(\log\frac{\la_j}{|V|^2+1})_k=\Ga_{jk}^j=(\log\mathbf{a}_j)_k,\
j\neq k$ we get $\frac{\la_j}{|V|^2+1}=\phi_j(u^j)\mathbf{a}_j$;
after a change of the $u^j$ variable into itself we can make
$\frac{4\la_j}{|V|^2+1}=\mathbf{a}_j,\ j=1,...,n$, so from
(\ref{eq:syst0}) $|\La|^2=-H(|V|^2+1)^2$ and
$\La^Td\La=-2dV^T(I_{1,n}+Ve_{n+1}^T)A(2V+(|V|^2-1)e_{n+1})=-2\La^T\del
R^T(I_{1,n}+Ve_{n+1}^T)A(2V+(|V|^2-1)e_{n+1})$. Thus again with
$\om:=\sum_{j=1}^n(e_je_j^TR^TR_j\del+\del R^TR_je_je_j^T)=-\om^T$
we have $d\La=\om\La-2\del
R^T(I_{1,n}+Ve_{n+1}^T)A(2V+(|V|^2-1)e_{n+1})$ and imposing
$d\wedge d\La=0$ we obtain the differential system
\begin{eqnarray}\label{eq:defqc}
d\wedge\om-\om\wedge\om=-4\del
R^T(I_{1,n}+Ve_{n+1}^T)A(I_{1,n}+e_{n+1}V^T)R\wedge\del,\
\om\wedge\del-\del\wedge R^TdR=0,\nonumber\\ dV=R\del\La,\
d\La=\om\La-2\del R^T(I_{1,n}+Ve_{n+1}^T)A(2V+(|V|^2-1)e_{n+1}),\nonumber\\
\La^T\La=-(2V+(|V|^2-1)e_{n+1})^TA(2V+(|V|^2-1)e_{n+1}),\nonumber\\
\om:=\sum_{j=1}^n(e_je_j^TR^TR_j\del+\del R^TR_je_je_j^T),\
R\subset\mathbf{O}_n(\mathbb{C})
\end{eqnarray}
in involution (that is no further conditions appear if one imposes
$d\wedge$ conditions and one uses the equations of the system
itself).

\section{The B\"{a}cklund transformation}

\subsection{(Isotropic) quadrics without center}\noindent

\noindent

Recall that the Ivory affinity between confocal (I)QWC is given by
$x_z=\sqrt{R_z}x_0+C(z)=\sqrt{R_z}LZ+C(z)=L(I_{1,n}\sqrt{R'_z}Z+
e_{n+1}(-I_{1,n}L^{-1}C(z)+e_{n+1})^TZ+L^{-1}C(z)),\
R_z:=I_{n+1}-zA,\
C(z):=-(\frac{1}{2}\int_0^z(\sqrt{R_w})^{-1}dw)B, A$ SJ,
$\ker(A)=\mathbb{C}e_{n+1},\ B=-e_{n+1}$ for QWC,
$\ker(A)=\mathbb{C}f_1,\ B=-\bar f_1$ for IQWC. Also
$L:=(\sqrt{A+e_{n+1}e_{n+1}^T})^{-1}$ for QWC and for IQWC $L$
satisfies $Le_{n+1}=f_1,\ L^T(A+\bar f_1\bar f_1^T)L=I_{n+1},\
A':=L^TA^2L=I_{1,n}(L^TL)^{-1}I_{1,n}$ SJ.

Consider two points $x_0^0,\ x_0^1\in x_0$ such that $x_0^0,\
x_z^1$ are in the symmetric TC
\begin{eqnarray}\label{eq:tcqwc}
0=(N_0^0)^T(x_z^1-x_0^0)\Leftrightarrow\
V_1^T\sqrt{R'_z}V_0-\frac{|V_1|^2+|V_0|^2}{2}+(V_0+V_1)^TL^{-1}C(z)-e_{n+1}^TL^{-1}C(z)=0
\Leftrightarrow\nonumber\\
|\sqrt{R'_z}V_1-V_0+I_{1,n}L^{-1}C(z)|^2=-zH_1 \Leftrightarrow\
x_z^1=x_0^0+[x_{0v_0^1}^0\ \ ...\ \
x_{0v_0^n}^0](\sqrt{R'_z}V_1-V_0+I_{1,n}L^{-1}C(z)).\nonumber\\
\end{eqnarray}
All these algebraic identities (and other needed later) boil down
to those established when a parametrization on (I)QWC was
introduced. Thus among the $2n$ functionally independent variables
$\{v_0^j,v_1^j\}_{j=1,...,n}$ a quadratic functional relation is
established and only $2n-1$ among them remain functionally
independent:
\begin{eqnarray}\label{eq:dtcqwc}
dV_0^T(\sqrt{R'_z}V_1-V_0+I_{1,n}L^{-1}C(z))=-dV_1^T(\sqrt{R'_z}V_0-V_1+I_{1,n}L^{-1}C(z)).
\end{eqnarray}
Given a deformation $x^0\subset\mathbb{C}^{2n-1}$ of $x_0^0$ (that
is $|dx^0|^2=|dx_0^0|^2$) with orthonormal normal frame
$N^0:=[N_{n+1}^0\ \ ...\ \ N_{2n-1}^0]$ consider the
$n$-dimensional sub-manifold
\begin{eqnarray}\label{eq:x1qwc}
x^1=x^0+[x_{v_0^1}^0\ \ ...\ \
x_{v_0^n}^0](\sqrt{R'_z}V_1-V_0+I_{1,n}L^{-1}C(z))\subset\mathbb{C}^{2n-1}
\end{eqnarray}
(that is we restrict $\{v_1^j\}_{j=1,...,n}$ to depend only on the
functionally independent $\{v_0^j\}_{j=1,...,n}$ and constants in
a manner that will subsequently become clear when we shall impose
the ACPIA).

Using $dx_{v_0^j}^0=\mathcal{V}^0dv_0^j-N^0(dN^0)^Tx_{v_0^j}^0$
and (\ref{eq:dtcqwc}) we have
\begin{eqnarray}\label{eq:dx1qwc}
dx_z^1=[-\mathcal{V}_0^0(\sqrt{R'_z}V_0-V_1+I_{1,n}L^{-1}C(z))^T+[x_{0v_0^1}^0\
\ ...\ \
x_{0v_0^n}^0]\sqrt{R'_z}]dV_1-N_0^0(dN_0^0)^T(x_z^1-x_0^0),\nonumber\\
dx^1=[-\mathcal{V}^0(\sqrt{R'_z}V_0-V_1+I_{1,n}L^{-1}C(z))^T+[x_{v_0^1}^0\
\ ...\ \ x_{v_0^n}^0]\sqrt{R'_z}]dV_1-
N^0(dN^0)^T(x^1-x^0).\nonumber\\
\end{eqnarray}
Since our intent is for $x^1$ to be a B transform (leaf) of the
seed $x^0$ we impose the ACPIA $|dx^1|^2=|dx_0^1|^2$. Keeping in
mind $|dx_z|^2=|dx_0|^2-z|dV|^2$ we obtain
$$z|dV_1|^2=|\begin{bmatrix}-i(dN_0^0)^T(x_z^1-x_0^0)\\
-(dN^0)^T(x^1-x^0)\end{bmatrix}|^2.$$
We take advantage now of the
conjugate system $(u^1,...,u^n)$ common to $x_0^0,x^0$ and of the
non-degenerate joined second fundamental forms property; according
to the principle of symmetry $0\leftrightarrow 1$ we would like
$(u^1,...,u^n)$ to be conjugate system to both $x^1$ and $x_0^1$
and also that the non-degenerate joined second fundamental forms
property holds for $x_0^1,x^1$.

We augment the column vector $-(dN^0)^T(x^1-x^0)$ with
$-i(dN_0^0)^T(x_z^1-x_0^0)$ as the first entry to obtain
$\begin{bmatrix}-i(dN_0^0)^T(x_z^1-x_0^0)\\-(dN^0)^T(x^1-x^0)\end{bmatrix}
=\begin{bmatrix}i(N_0^0)^T[dx_{0v_0^1}^0\ \ ...\ \
dx_{0v_0^n}^0]\\
(N^0)^T[dx_{v_0^1}^0\ \ ...\ \
dx_{v_0^n}^0]\end{bmatrix}(\sqrt{R'_z}V_1-V_0+I_{1,n}L^{-1}C(z))
=\begin{bmatrix}i(N_0^0)^T[dx_{0u^1}^0\ \ ...\ \
dx_{0u^n}^0]\\
(N^0)^T[dx_{u^1}^0\ \ ...\ \
dx_{u^n}^0]\end{bmatrix}\frac{\pa(u^1,...,u^n)}{\pa(v_0^1,...,v_0^n)}
(\sqrt{R'_z}V_1-V_0+I_{1,n}L^{-1}C(z))
=S_0\mathrm{diag}[\la_{01}du^1\ \ ...\ \ \la_{0n}du^n]\\
\frac{\pa(u^1,...,u^n)}{\pa(v_0^1,...,v_0^n)}(\sqrt{R'_z}V_1-V_0+I_{1,n}L^{-1}C(z))=S_0\del
R_0^T(\sqrt{R'_z}V_1-V_0+I_{1,n}L^{-1}C(z))$, where
$S_0\subset\mathbf{O}_n(\mathbb{C})$ is the orthogonal matrix with
columns $\la_{0k}^{-1}[ih_k^0\ \ h_k^{n+1}\ \ ...\ \
h_k^{2n-1}]^T,\ \{h_k^0\}_k,\ \{h_k^{\al}\}_{k,\al}$ are the
second fundamental forms of $x_0^0,x^0$ and
$\frac{\pa(v_0^1,...,v_0^n)}{\pa(u^1,...,u^n)}=R_0\La_0$. Thus
$dV_1=-\frac{1}{\sqrt{z}}M_1S_0\del
R_0^T(\sqrt{R'_z}V_1-V_0+I_{1,n}L^{-1}C(z))=R_1\del\La_1,\
M_1\subset\mathbf{O}_n(\mathbb{C}),\
R_1:=M_1S_0\subset\mathbf{O}_n(\mathbb{C}),\
\La_1:=-\frac{1}{\sqrt{z}}R_0^T(\sqrt{R'_z}V_1-V_0+I_{1,n}L^{-1}C(z))$.
For the prime integral property $|\La_1|^2=-H_1$ we have
$|\sqrt{R'_z}V_1-V_0+I_{1,n}L^{-1}C(z)|^2=-zH_1$. Note that $S_0$
is undetermined up to multiplication on the left with
$R'\subset\mathbf{O}_n(\mathbb{C}),\ R'e_1=e_1$, so the same
statement on the right holds for $M_1$, but $M_1S_0$ is well
defined.

Assuming that $V_1,\ \La_1$ can be found with the properties
needed so far, (\ref{eq:syst0}) is satisfied and the existence of
the B transformation is proved, provided we exhibit an explicit
algebraic transformation of the solution $V_0,\ \La_0$ of
(\ref{eq:systqwc}) to the solution $V_1,\ \La_1$ of the same
system, with the transformation matrix depending only on the
orthogonal $R_0,R_1$ solutions of (\ref{eq:defqwc}): this
dependence will reveal the analytic B transformation between
solutions $R_0,R_1$ of (\ref{eq:defqwc}).

We have
\begin{eqnarray}\label{eq:rlaqwc}
\sqrt{z}R_1\La_0=\sqrt{R'_z}V_0-V_1+I_{1,n}L^{-1}C(z),\nonumber\\
-\sqrt{z}R_0\La_1=\sqrt{R'_z}V_1-V_0+I_{1,n}L^{-1}C(z)
\end{eqnarray}
(the symmetry $(0,\sqrt{z})\leftrightarrow(1,-\sqrt{z})$ is
required by (\ref{eq:dtcqwc})). From the first relation of
(\ref{eq:rlaqwc}) we get $V_1$ as an algebraic expression of
$R_1,\ V_0,\ \La_0$; replacing this into the second relation of
(\ref{eq:rlaqwc}) we get $\La_1$ as an algebraic expression of
$R_0,\ R_1,\ V_0,\
\La_0$:
\begin{eqnarray}\label{eq:algqwc}
V_1=\sqrt{R'_z}V_0+I_{1,n}L^{-1}C(z)-\sqrt{z}R_1\La_0,\nonumber\\
\La_1=R_0^T(\sqrt{z}A'V_0+\sqrt{R'_z}R_1\La_0+\sqrt{z}I_{1,n}L^{-1}B).
\end{eqnarray}
Differentiating the first relation of (\ref{eq:algqwc}) and using
the second we get $R_1\del
R_0^T(\sqrt{z}A'V_0+\sqrt{R'_z}R_1\La_0+\sqrt{z}I_{1,n}L^{-1}B)=\sqrt{R'_z}R_0\del\La_0
-\sqrt{z}[dR_1\La_0+R_1(\om_0\La_0-\del R_0^T(A'V_0+L^{-1}B))]$,
or the Ricatti equation
\begin{eqnarray}\label{eq:bqwc}
-dR_1=R_1\om_0+R_1\del R_0^TDR_1-DR_0\del,\
 D:=\frac{\sqrt{R'_z}}{\sqrt{z}}
\end{eqnarray}
in $R_1$.

Again this is similar to Tenenblat-Terng's B transformation in
Terng's interpretation. Imposing the compatibility condition
$d\wedge$ on (\ref{eq:bqwc}) and using the equation itself we thus
need: $0=-(R_1\om_0+R_1\del R_0^TDR_1-DR_0\del)\wedge\om_0
+R_1d\wedge\om_0-(R_1\om_0+ R_1\del R_0^TDR_1-DR_0\del)\wedge\del
R_0^TDR_1-R_1\del\wedge dR_0^TDR_1+R_1\del
R_0^TD\wedge(R_1\om_0+R_1\del R_0^TDR_1-DR_0\del)
-DdR_0\wedge\del=R_1(-\om_0\wedge\om_0 +d\wedge\om_0+\del
R_0^TA'R_0\wedge\del)+ DR_0(\del\wedge\om_0-R_0^TdR_0\wedge\del)
-R_1(\om_0\wedge\del-\del\wedge R_0^TdR_0)R_0^TDR_1=0$ because
$R_0$ satisfies (\ref{eq:defqwc}). Therefore (\ref{eq:bqwc}) is
completely integrable and admits solution for any initial value of
$R_1$. If the initial value is orthogonal, we would like the
solution to remain orthogonal:
$d(R_1R_1^T-I_n)=dR_1R_1^T+R_1dR_1^T=-(R_1\om_0+R_1\del R_0^TDR_1
-DR_0\del)R_1^T-R_1(-\om_0R_1^T+R_1^TDR_0\del R_1^T-\del R_0^TD)=
-R_1\del R_0^TD(R_1R_1^T-I_n)- (R_1R_1^T-I_n)DR_0\del R_1^T$, so
$R_1R_1^T-I_n$ is a solution of a linear differential equation and
remains $0$ if initially it was $0$. The fact that $R_1$ is itself
a solution of (\ref{eq:defqwc}) follows from the symmetry
$(0,\sqrt{z})\leftrightarrow(1,-\sqrt{z})$ and the fact that
$d\wedge dR_0=0$ (basically we use the converse of the proven
results).

Summing up:

{\it Given $R_0\subset\mathbf{O}_n(\mathbb{C})$ solution of
(\ref{eq:defqwc}) and $R_1\subset\mathbf{M}_n(\mathbb{C})$
solution of the Ricatti equation (\ref{eq:bqwc}), then $R_1$
remains orthogonal if initially it was orthogonal and in this case
it is another solution of (\ref{eq:defqwc}) (thus producing an
$(\frac{n(n-1)}{2}+1)$-dimensional family of solutions). Moreover
if $V_0,\ \La_0$ are solutions of (\ref{eq:systqwc}) associated to
$R_0$ (thus producing a seed deformation
$x^0\subset\mathbb{C}^{2n-1}$ of $x_0^0$), then $V_1,\ \La_1$
given by
\begin{eqnarray}\label{eq:algsystqwc}
\begin{bmatrix}V_1\\\La_1\end{bmatrix}=\sqrt{z}
\begin{bmatrix}I_n&0\\0&R_0^T\end{bmatrix}
(\begin{bmatrix}D&-I_n\\A'&D\end{bmatrix}
\begin{bmatrix}I_n&0\\0&R_1\end{bmatrix}
\begin{bmatrix}V_0\\\La_0\end{bmatrix}+
\begin{bmatrix}I_{1,n}\frac{L^{-1}C(z)}{\sqrt{z}}\\
I_{1,n}L^{-1}B\end{bmatrix}),\ D:=\frac{\sqrt{R'_z}}{\sqrt{z}}
\end{eqnarray}
are solutions of (\ref{eq:systqwc}) associated to $R_1$ (thus
producing a leaf deformation $x^1\subset\mathbb{C}^{2n-1}$ of
$x_0^1$) and we have the symmetry
$(0,\sqrt{z})\leftrightarrow(1,-\sqrt{z})$.}

Again one can easily develop the algebraic formula of the BPT at
the level of algebraic transformation (\ref{eq:algsystqwc}) of
solutions (\ref{eq:bqwc}) and then prove that indeed it induces
the BPT of solutions of (\ref{eq:defqwc}):
$\begin{bmatrix}I_n&0\\0&R_2R_1^T\end{bmatrix}
(\begin{bmatrix}D_2&-R_3R_0^T\\A'&D_2R_3R_0^T\end{bmatrix}
(\begin{bmatrix}D_1&-I_n\\A'&D_1\end{bmatrix}
\begin{bmatrix}V_0\\R_1\La_0\end{bmatrix}+
\begin{bmatrix}I_{1,n}\frac{L^{-1}C(z_1)}{\sqrt{z_1}}\\
I_{1,n}L^{-1}B\end{bmatrix})+\frac{1}{\sqrt{z_1}}
\begin{bmatrix}I_{1,n}\frac{L^{-1}C(z_2)}{\sqrt{z_2}}\\
I_{1,n}L^{-1}B\end{bmatrix})=
\begin{bmatrix}D_1&-R_3R_0^T\\A'&D_1R_3R_0^T\end{bmatrix}
(\begin{bmatrix}D_2&-R_2R_1^T\\A'&D_2R_2R_1^T\end{bmatrix}
\begin{bmatrix}V_0\\R_1\La_0\end{bmatrix}+
\begin{bmatrix}I_{1,n}\frac{L^{-1}C(z_2)}{\sqrt{z_2}}\\
I_{1,n}L^{-1}B\end{bmatrix})+\frac{1}{\sqrt{z_2}}
\begin{bmatrix}I_{1,n}\frac{L^{-1}C(z_1)}{\sqrt{z_1}}\\
I_{1,n}L^{-1}B\end{bmatrix}$, or
\begin{eqnarray}\label{eq:bptqwc}
R_3R_0^T=(D_2-D_1R_2R_1^T) (D_2R_2R_1^T-D_1)^{-1}.
\end{eqnarray}
The formula thus established the remaining part will follow in a
manner similar to Terng's approach. First we need
$R_3R_0^T\subset\mathbf{O}_n(\mathbb{C})$, which follows from
$(D_2-R_2R_1^TD_1) (D_2R_2R_1^T-D_1)= (R_2R_1^TD_2-D_1)
(D_2-D_1R_2R_1^T)$.

Therefore we only need to prove that $R_3$ given by
(\ref{eq:bptqwc}) satisfies (\ref{eq:bqwc}) for $(R_0,z)$ replaced
by $(R_1,z_2),\ (R_2,z_1)$; by symmetry it is enough to prove only
one relation. Since $-dR_1=R_1\om_0+R_1\del
R_0^TD_1R_1-D_1R_0\del,\ -dR_2=R_2\om_0+R_2\del
R_0^TD_2R_2-D_2R_0\del$, we get $d(R_2R_1^T)=-(R_2\om_0+R_2\del
R_0^TD_2R_2-D_2R_0\del)R_1^T -R_2(-\om_0R_1^T+R_1^TD_1R_0\del
R_1^T-\del R_0^TD_1)=-(R_2R_1^T)R_1\del R_0^T(D_2R_2R_1^T-D_1)
+(D_2-R_2R_1^TD_1)R_0\del R_1^T$. Thus if we prove the similar
relation $d(R_3R_0^T)=\\-(R_3R_0^T)R_0\del R_1^T(D_2R_3R_0^T+D_1)
+(D_2+R_3R_0^TD_1)R_1\del R_0^T$, then since
$-dR_0=R_0\om_1-R_0\del R_1^TD_1R_0+D_1R_1\del$ we obtain what we
want: $-dR_3=R_3\om_1+R_3\del R_1^TD_2R_3-D_2R_1\del$.
Differentiating (\ref{eq:bptqwc}) we get
$d(R_3R_0^T)(D_2R_2R_1^T-D_1)= -(R_3R_0^TD_2+D_1)d(R_2R_1^T)$;
thus we need to prove $\\\ [-(R_3R_0^T)R_0\del
R_1^T(D_2R_3R_0^T+D_1) +(D_2+R_3R_0^TD_1)R_1\del
R_0^T](D_2R_2R_1^T-D_1)= -(R_3R_0^TD_2+D_1)[-(R_2R_1^T)R_1\del
R_0^T(D_2R_2R_1^T-D_1) +(D_2-R_2R_1^TD_1)R_0\del R_1^T]$. The
terms containing $R_1\del R_0^T$ become $D_2+R_3R_0^TD_1=
(R_3R_0^TD_2+D_1)R_2R_1^T$ which follows directly from
(\ref{eq:bptqwc}) and the terms containing $R_0\del R_1^T$ become
$R_0\del R_1^T(D_2R_3R_0^T+D_1) (D_2R_2R_1^T-D_1)=(R_3R_0^T)^T
(R_3R_0^TD_2+D_1) (D_2-R_2R_1^TD_1)R_0\del R_1^T$ which follows
from $(D_2R_3R_0^T+D_1) (D_2R_2R_1^T-D_1)=D_2^2-D_1^2=
(\frac{1}{z_2}-\frac{1}{z_1})I_n=(R_3R_0^T)^T(R_3R_0^TD_2+D_1)
(D_2-R_2R_1^TD_1)$.

To prove the existence of the $\mathcal{M}_p$ configuration we
need only prove the existence of the $\mathcal{M}_3$ configuration
(discrete deformations in $\mathbb{C}^{2n-1}$ of $x_0$ will be
obtained by considering $\mathbb{Z}^n$ lattices of B
transformations which are subsets of infinite M\"{o}bius
configurations, similarly to Bobenko-Pinkall's approach
\cite{BP}).

Consider $(D_1D_2)^{-1}[(D_1^2-D_2^2)
D_1R_1(D_1R_1-D_2R_2)^{-1}-D_1^2]=R_3R_0^{-1}=(D_1D_2)^{-1}[(D_1^2-D_2^2)
D_2R_2(D_1R_1-D_2R_2)^{-1}-D_2^2],\
R_5R_0^{-1}=(D_1D_3)^{-1}[(D_3^2-D_1^2)
D_1R_1(D_3R_4-D_1R_1)^{-1}-D_1^2],\
R_6R_0^{-1}=(D_2D_3)^{-1}[(D_2^2-D_3^2)
D_2R_2(D_2R_2-D_3R_4)^{-1}-D_2^2]$; thus with
$\Box:=(D_2^2-D_3^2)D_1R_1+(D_3^2-D_1^2)D_2R_2+(D_1^2-D_2^2)D_3R_4$
we have $(D_2R_3R_0^{-1}-D_3R_5R_0^{-1})^{-1}R_1=[(D_1^2-D_2^2)
(D_1R_1-D_2R_2)^{-1}-(D_3^2-D_1^2)(D_3R_4-D_1R_1)^{-1}]^{-1}
=(D_1R_1-D_2R_2)\Box^{-1}(D_3R_4-D_1R_1)$ and similarly
$(D_3R_6R_0^{-1}-D_1R_3R_0^{-1})^{-1}R_2=(D_1R_1-D_2R_2)\Box^{-1}(D_2R_2-D_3R_4)$.
Now $D_1[(D_2^2-D_3^2)D_2R_3R_0^{-1}\\
(D_2R_3R_0^{-1}-D_3R_5R_0^{-1})^{-1}R_1-D_2^2R_1]=
(D_1^2D_2R_2-D_2^2D_1R_1)\Box^{-1}(D_2^2-D_3^2)(D_3R_4-D_1R_1)-D_2^2D_1R_1
=(D_1^2D_2R_2-D_2^2D_1R_1)\Box^{-1}(D_3^2-D_1^2)(D_2R_2-D_3R_4)-D_1^2D_2R_2
=D_2[(D_3^2-D_1^2)D_1R_3R_0^{-1}\\
(D_3R_6R_0^{-1}-D_1R_3R_0^{-1})^{-1}R_2-D_1^2R_2]$, so the very
lhs and rhs provide the good definition of and afford themselves
the name $D_1D_2D_3R_7$.

We shall now show that the facets $dx^1$, rigidly transported by
rolling to facets centered at $x_z^1$ when $x^0$ rolls on $x_0^0$,
will cut the tangent spaces of $x_z^1$ along rulings.

In this case it is useful to extend the orthonormal normal frame
$N^0$ with a $0$ column vector on the left and call the extended
'orthonormal normal frame' thus obtained still orthonormal frame;
all computations will remain valid as long as one uses transpose
instead of inverse. The facet $dx^1$ will be transported to the
facet $dx_z^1+N_0^0(dN_0^0)^T(x_z^1-x_0^0)+\sqrt{z}[0\ N_0^0\ \
e_{n+2}\ \ ...\ \ e_{2n-1}]MdV_1=([x_{zv_1^1}^1\ \ ...\ \
x_{zv_1^n}^1]+\frac{N_0^0}{\sqrt{H_0}}(\sqrt{R'_z}V_0-V_1+I_{1,n}L^{-1}C(z))^T+\sqrt{z}[0\
N_0^0\ \ e_{n+2}\ \ ...\ \ e_{2n-1}]M)dV_1,\
dV_1=\frac{1}{\sqrt{z}}R_1\del
R_0^T(\sqrt{R'_z}V_1-V_0+I_{1,n}L^{-1}C(z))$, where
$M\subset\mathbf{O}_n(\mathbb{C})$ is determined up to
multiplication on the left with an orthogonal matrix which fixes
$e_1$ (due to the indeterminacy of the rolling in the normal
bundle) by $M^Te_1=\pm
i\frac{\sqrt{R'_z}V_0-V_1+I_{1,n}L^{-1}C(z)}{\sqrt{zH_0}}$
(because of the ACPIA requirement
$z(e_1^TMdV_1)^2=-[(dN_0^0)^T(x_z^1-x_0^0)]^2=
-\frac{1}{H_0}((\sqrt{R'_z}V_0-V_1+I_{1,n}L^{-1}C(z))^TdV_1)^2$).
The vectors of the facet are obtained by replacing $\del$ with
diagonal constant matrices $\Del$; let
$\Del':=\frac{1}{\sqrt{z}}R_1\Del
R_0^T(\sqrt{R'_z}V_1-V_0+I_{1,n}L^{-1}C(z))$. Thus we need
$e_j^TM\Del'=0,\ j=3,...,n$, so
$M\Del'=e_1^TM\Del'e_1+e_2^TM\Del'e_2$; from $(N_z^1)^Tw=0$ we
obtain $e_2^TM\Del'=\pm ie_1^TM\Del'$; after scaling we obtain
$w=\pm[x_{zv_1^1}^1\ \ ...\ \ x_{zv_1^n}^1](e_1\pm iM^Te_2)$. Now
the condition $w^TAR_z^{-1}w=0$ that $w$ is a ruling on $x_z^1$
becomes $|e_1\pm iM^Te_2|^2=0$ and it is thus satisfied.

Now the notion of W congruence for $n=2$ requires that the
asymptotic directions on the two focal surfaces correspond, so the
natural generalization for $n\ge 3$ is that the asymptotic
directions on the two focal sub-manifolds $x^0,x^1$ should
correspond. An asymptotic direction $v=a^jx_j$ must satisfy
$a^ja^kh^{\al}_{jk}=0,\ \al=n+1,...,2n-1$; in our case
(restricting to conjugate system) we need
$\sum_j(a^j)^2\la_j\frac{h^{\al}_j}{\la_j}=0$; keeping account of
the orthogonal properties of the non-degenerate joined second
fundamental forms we get $a^j=\pm 1$ and thus they are the same
for $x^0,x^1$.

\subsection{Quadrics with center}\noindent

\noindent Recall that we have the Ivory affinity
$x_z=\sqrt{R_z}(\sqrt{A})^{-1}X,\
X=\frac{2V+(|V|^2-1)e_{n+1}}{|V|^2+1}$. Note
$$\frac{(|V|^2+1)^2}{4}\sum_{k=1}^nx_{zv^k}x_{zv^k}^T+x_zx_z^T=A^{-1}-zI_{n+1}.$$

Consider two points $x_0^0,\ x_0^1\in x_0$ such that $x_0^0,\
x_z^1$ are in the symmetric TC
\begin{eqnarray}\label{eq:tcqc}
(x_z^1-x_0^0)^TN_0^0=0\Leftrightarrow
X_1^T\sqrt{R_z}X_0=1\Leftrightarrow
x_z^1=x_0^0+\frac{(|V_0|^2+1)^2}{4}\sum_{k=1}^nX_1^T\sqrt{R_z}X_{0v_0^k}x_{0v_0^k}^0
\Leftrightarrow\nonumber\\
(2V_0+(|V_0|^2-1)e_{n+1})^T\sqrt{R_z}(2V_1+(|V_1|^2-1)e_{n+1})=(|V_0|^2+1)(|V_1|^2+1).
\end{eqnarray}
Thus among the $2n$ functionally independent variables
$\{v_0^j,v_1^j\}_{j=1,...,n}$ a quadratic functional relation is
established and only $2n-1$ among them remain functionally
independent:
\begin{eqnarray}\label{eq:dtcqc}
X_1^T\sqrt{R_z}dX_0=-X_0^T\sqrt{R_z}dX_1.
\end{eqnarray}
Given a deformation $x^0\subset\mathbb{C}^{2n-1}$ of $x_0^0$ (that
is $|dx^0|^2=|dx_0^0|^2$) with orthonormal normal frame
$N^0:=[N_{n+1}^0\ \ ...\ \ N_{2n-1}^0]$ consider the
$n$-dimensional sub-manifold
\begin{eqnarray}\label{eq:x1qc}
x^1=x^0+\frac{(|V_0|^2+1)^2}{4}\sum_{k=1}^nX_1^T\sqrt{R_z}X_{0v_0^k}x_{v_0^k}^0
\subset\mathbb{C}^{2n-1}
\end{eqnarray}
(that is we restrict $\{v_1^j\}_{j=1,...,n}$ to depend only on the
functionally independent $\{v_0^j\}_{j=1,...,n}$ and constants in
a manner that will subsequently become clear when we shall impose
the ACPIA).

Using $dx_{v_0^k}^0=-\frac{\pa\log(|V_0|^2+1)}{\pa
v_0^k}dx^0-d\log(|V_0|^2+1)x_{v_0^k}^0+(\mathcal{V}^0+\sum_{j=1}^n\frac{\pa\log(|V_0|^2+1)}{\pa
v_0^j}x_{v_0^j}^0)dv_0^k-N^0(dN^0)^Tx_{v_0^k}^0$ and
(\ref{eq:dtcqc}) we have
\begin{eqnarray}\label{eq:dx1qc}
dx_z^1=(-\mathcal{V}_0^0X_0^T+\frac{(|V_0|^2+1)^2}{4}\sum_{k=1}^nx_{0v_0^k}^0X_{0v_0^k}^T)
\sqrt{R_z}dX_1-N_0^0(dN_0^0)^T(x_z^1-x_0^0),\nonumber\\
dx^1=(-\mathcal{V}^0X_0^T+\frac{(|V_0|^2+1)^2}{4}\sum_{k=1}^nx_{v_0^k}^0X_{0v_0^k}^T)
\sqrt{R_z}dX_1-N^0(dN^0)^T(x^1-x^0).
\end{eqnarray}
Since our intent is for $x^1$ to be a B transform (leaf) of the
seed $x^0$ we impose the ACPIA $|dx^1|^2=|dx_0^1|^2$. Keeping in
mind $|dx_z|^2=|dx_0|^2-z|dX|^2$ we obtain
$$z|dX_1|^2=|\begin{bmatrix}-i(dN_0^0)^T(x_z^1-x_0^0)\\-(dN^0)^T(x^1-x^0)\end{bmatrix}|^2.$$

We take advantage now of the conjugate system $(u^1,...,u^n)$
common to $x_0^0,x^0$ and of the non-degenerate joined second
fundamental forms property; according to the principle of symmetry
$0\leftrightarrow 1$ we would like $(u^1,...,u^n)$ to be conjugate
system to both $x^1$ and $x_0^1$ and also that the non-degenerate
joined second fundamental forms property holds for $x_0^1,x^1$. We
obtain by computations similar to the (I)QWC case
$dV_1=R_1\del\La_1,\
\La_1:=\frac{(|V_1|^2+1)(|V_0|^2+1)}{2\sqrt{z}}R_0^T[X_1^T\sqrt{R_z}X_{0v_0^1}\
\ ...\ \ X_1^T\sqrt{R_z}X_{0v_0^n}]^T$ (note that the prime
integral property $|\La_1|^2=-H_1(|V_1|^2+1)^2$ is satisfied).

Thus
\begin{eqnarray}\label{eq:rlaqc}
\sqrt{z}R_0\La_1=(I_{1,n}+V_0e_{n+1}^T)\sqrt{R_z}(2V_1+(|V_1|^2-1)e_{n+1})
-V_0(|V_1|^2+1),\nonumber\\
-\sqrt{z}R_1\La_0=(I_{1,n}+V_1e_{n+1}^T)\sqrt{R_z}(2V_0+(|V_0|^2-1)e_{n+1})
-V_1(|V_0|^2+1) \end{eqnarray} (the symmetry
$(0,\sqrt{z})\leftrightarrow(1,-\sqrt{z})$ follows from
(\ref{eq:dtcqc})). From the second relation of (\ref{eq:rlaqc}) we
get $V_1$ as an algebraic expression of $R_1,\ V_0,\ \La_0$;
replacing this into the first relation of (\ref{eq:rlaqc}) we get
$\La_1$ as an algebraic expression of $R_0,\ R_1,\ V_0,\ \La_0$:
\begin{eqnarray}\label{eq:algqc}
V_1=-\frac{\sqrt{z}R_1\La_0+I_{1,n}\sqrt{R_z}(2V_0+(|V_0|^2-1)e_{n+1})}
{e_{n+1}^T\sqrt{R_z}(2V_0+(|V_0|^2-1)e_{n+1})-|V_0|^2-1},\nonumber\\
\La_1=2R_0^T\frac{(I_{1,n}+V_0e_{n+1}^T)[\sqrt{z}A(2V_0+(|V_0|^2-1)e_{n+1})
-\sqrt{R_z}(I_{1,n}+e_{n+1}V_1^T)R_1\La_0]+V_0V_1^TR_1\La_0}
{e_{n+1}^T\sqrt{R_z}(2V_0+(|V_0|^2-1)e_{n+1})-|V_0|^2-1}.\nonumber\\
\end{eqnarray}
Differentiating the second relation of (\ref{eq:rlaqc}) we get
$-\sqrt{z}dR_1\La_0-\sqrt{z}R_1[\om_0\La_0-2\del
R_0^T(I_{1,n}+V_0e_{n+1}^T)A(2V_0\\+(|V_0|^2-1)e_{n+1})]-
2I_{1,n}\sqrt{R_z}(I_{1,n}+e_{n+1}V_0^T)R_0\del\La_0=
[e_{n+1}^T\sqrt{R_z}(2V_0+(|V_0|^2-1)e_{n+1})-|V_0|^2-1]R_1\del\La_1
+2V_1[e_{n+1}^T\sqrt{R_z}(I_{1,n}+e_{n+1}V_0^T)-V_0^T]R_0\del\La_0$,
or equivalently the Ricatti equation
\begin{eqnarray}\label{eq:bqc}
-dR_1=R_1\om_0+2I_{1,n}\frac{\sqrt{R_z}}{\sqrt{z}}(I_{1,n}+e_{n+1}V_0^T)R_0\del
-2R_1\del
R_0^T(I_{1,n}+V_0e_{n+1}^T)\frac{\sqrt{R_z}}{\sqrt{z}}R_1\nonumber\\
+2R_1\del R_0^T\frac{[(I_{1,n}+V_0e_{n+1}^T)\sqrt{R_z}e_{n+1}-V_0]
[\La_0^T+(2V_0^T+(|V_0|^2-1)e_{n+1}^T)\frac{\sqrt{R_z}}{\sqrt{z}}R_1]}
{e_{n+1}^T\sqrt{R_z}(2V_0+(|V_0|^2-1)e_{n+1})-|V_0|^2-1}\nonumber\\
-2\frac{[R_1\La_0+I_{1,n}\frac{\sqrt{R_z}}{\sqrt{z}}(2V_0+(|V_0|^2-1)e_{n+1})]
[e_{n+1}^T\sqrt{R_z}(I_{1,n}+e_{n+1}V_0^T)-V_0^T]}
{e_{n+1}^T\sqrt{R_z}(2V_0+(|V_0|^2-1)e_{n+1})-|V_0|^2-1}R_0\del
\end{eqnarray}
in $R_1$.

Note that with
$M:=I_{1,n}\frac{\sqrt{R_z}}{\sqrt{z}}(I_{1,n}+e_{n+1}V_0^T),\
N:=I_{1,n}\frac{\sqrt{R_z}}{\sqrt{z}}(2V_0+(|V_0|^2-1)e_{n+1}),\\
W:=(I_{1,n}+V_0e_{n+1}^T)\sqrt{R_z}e_{n+1}-V_0,\
U:=e_{n+1}^T\sqrt{R_z}(2V_0+(|V_0|^2-1)e_{n+1})-|V_0|^2-1$ we have
$dN=2MdV_0,\ dU=2W^TdV_0,\ dM=CdV_0^T,\ dW=cdV_0$ and
(\ref{eq:bqc}) can be written as
$$-dR_1=R_1\om_0+2MR_0\del-2R_1\del R_0^TM^TR_1+\frac{2}{U}R_1\del
R_0^TW(\La_0^T+N^TR_1) -\frac{2}{U}(R_1\La_0+N)W^TR_0\del.$$
Imposing the $d\wedge$ condition on (\ref{eq:bqc}) and using the
equation itself we need $0=-[R_1\om_0+2MR_0\del-2R_1\del
R_0^TM^TR_1+\frac{2}{U}R_1\del R_0^TW(\La_0^T+N^TR_1)
-\frac{2}{U}(R_1\La_0+N)W^TR_0\del]\wedge[\om_0-2\del
R_0^TM^TR_1+\frac{2}{U}\del
R_0^TW(\La_0^T+N^TR_1)-\frac{2}{U}\La_0W^TR_0\del] -R_1\del
R_0^T(2M^T-\frac{2}{U}WN^T)\wedge[R_1\om_0+2MR_0\del-2R_1\del
R_0^TM^TR_1+\frac{2}{U}R_1\del R_0^TW(\La_0^T+N^TR_1)
-\frac{2}{U}(R_1\La_0+N)W^TR_0\del]+R_1d\wedge\om_0+2CdV_0^TR_0\wedge\del
+2MdR_0\wedge\del-2R_1\del\wedge
R_0^TdR_0R_0^TM^TR_1+2R_1\del\wedge R_0^TdV_0C^TR_1
+\frac{2}{U}R_1\del\wedge
R_0^TdR_0R_0^TW(\La_0^T+N^TR_1)-\frac{2c}{U}R_1\del\wedge
R_0^TdV_0(\La_0^T+N^TR_1) -\frac{2}{U}R_1\del\wedge
R_0^TW(d\La_0^T+2dV_0^TM^TR_1) +\frac{4}{U^2}R_1\del\wedge
R_0^TW(\La_0^T+N^TR_1)W^TdV_0
-\frac{2}{U}(R_1d\La_0+2MdV_0)W^TR_0\wedge\del-\frac{2}{U}(R_1\La_0+N)cdV_0^TR_0\wedge\del
+\frac{4}{U^2}(R_1\La_0+N)W^TW^TdV_0R_0\wedge\del-\frac{2}{U}(R_1\La_0+N)W^TdR_0\wedge\del=
\frac{4}{U^2}[W^T(dV_0-R_0\del\La_0)\wedge[(R_1\La_0+N)W^TR_0\del-R_1\del
R_0^TW(\La_0^T+N^TR_1)]+(|\La_0|^2-|N|^2)R_1\del
R_0^TWW^TR_0\wedge\del]+\frac{2}{U}[-cR_1\del
R_0^T\wedge(dV_0-R_0\del\La_0)(\La_0^T+N^TR_1)-c(R_1\La_0+N)(dV_0-R_0\del\La_0)^TR_0\wedge\del
-R_1\del\wedge
R_0^TW[\La_0^T\om_0+d\La_0^T+2(dV_0-R_0\del\La_0)^TM^TR_1]
+[R_1(\om_0\La_0-d\La_0)+2M(R_0\del\La_0-dV_0)]W^TR_0\wedge\del
+(R_1\La_0+N)W^TR_0(\del\wedge\om_0-R_0^TdR_0\wedge\del)
-R_1(\om_0\wedge\del-\del\wedge R_0^TdR_0)R_0^TW(\La_0^T+N^TR_1)
+2R_1\del R_0^T(WN^TM+M^TNW^T)R_0\wedge\del]
-2MR_0(\del\wedge\om_0-R_0^TdR_0\wedge\om_0)+R_1[-\om_0\wedge\om_0+d\wedge\om_0-4\del
R_0^TM^TMR_0\wedge\del+(\om_0\wedge\del-\del\wedge
R_0^TdR_0)R_0^TM^TR_1]+2C(dV_0-R_0\del\La_0)^TR_0\wedge\del+2R_1\del\wedge
R_0^T(dV_0-R_0\del\La_0)C^TR_1$. Using the fact that $V_0,\
\La_0,\ R_0$ are solutions of (\ref{eq:defqc}) we need $0=R_1\del
R_0^T[|\La_0|^2-|N|^2]WW^T
+UW[(2V_0+(|V_0|^2-1)e_{n+1})^TA(I_{1,n}+e_{n+1}V_0^T)+N^TM]
+U[(I_{1,n}+V_0e_{n+1}^T)A(2V_0+(|V_0|^2-1)e_{n+1})+M^TN]W^T
-U^2[M^TM+(I_{1,n}+V_0e_{n+1}^T)A(I_{1,n}+e_{n+1}V_0^T)]R_0\wedge\del$;
using
$|\La_0|^2-|N|^2=-\frac{1}{z}[(|V_0|^2+1)^2-(U+|V_0|^2+1)^2],\
(2V_0+(|V_0|^2-1)e_{n+1})^TA(I_{1,n}+e_{n+1}V_0^T)+N^TM=
-\frac{1}{z}[U(W+V_0)^T+(|V_0|^2+1)W^T],\
M^TM+(I_{1,n}+V_0e_{n+1}^T)A(I_{1,n}+e_{n+1}V_0^T=\frac{1}{z}[I_n-WW^T-WV_0^T-V_0W^T]$
this is straightforward. Therefore (\ref{eq:bqc}) is completely
integrable and admits solution for any initial value of $R_1$. If
the initial value is orthogonal, we would like the solution to
remain orthogonal:
$d(R_1R_1^T-I_n)=dR_1R_1^T+R_1dR_1^T=-[R_1\om_0+2MR_0\del-2R_1\del
R_0^TM^TR_1+\frac{2}{U}R_1\del R_0^TW(\La_0^T+N^TR_1)
-\frac{2}{U}(R_1\La_0+N)W^TR_0\del]R_1^T-R_1[-\om_0R_1^T+2\del
R_0^TM^T-2R_1^TMR_0\del R_1^T+\frac{2}{U}(\La_0+R_1^TN)W^TR_0\del
R_1^T -\frac{2}{U}\del
R_0^TW(\La_0^TR_1^T+N^T)]=2(R_1R_1^T-I_n)MR_0\del R_1^T+2R_1\del
R_0^TM^T(R_1R_1^T-I_n)-\frac{2}{U}R_1\del
R_0^TWN^T(R_1R_1^T-I_n)-\frac{2}{U}(R_1R_1^T-I_n)NW^TR_0\del
R_1^T$, so $R_1R_1^T-I_n$ is a solution of a linear differential
equation and remains $0$ if initially it was $0$. The fact that
$R_1$ is itself a solution of (\ref{eq:defqc}) follows from the
symmetry $(0,\sqrt{z})\leftrightarrow(1,-\sqrt{z})$ and the fact
that $d\wedge dR_0=0$ (basically we use the converse of the proven
results).

We shall now show that the facets $dx^1$, rigidly transported by
rolling to facets centered at $x_z^1$ when $x^0$ rolls on $x_0^0$,
will cut the tangent spaces of $x_z^1$ along rulings.

In this case it is useful to extend the orthonormal normal frame
$N^0$ with a $0$ column vector on the left and call the extended
'orthonormal normal frame' thus obtained still orthonormal frame;
all computations will remain valid as long as one uses transpose
instead of inverse. The facet $dx^1$ will be transported to the
facet $dx_z^1+N_0^0(dN_0^0)^T(x_z^1-x_0^0)+2\sqrt{z}[0\ N_0^0\ \
e_{n+2}\ \ ...\ \ e_{2n-1}]\frac{MdV_1}{|V_1|^2+1}=([x_{zv_1^1}^1\
\ ...\ \
x_{zv_1^n}^1]+\frac{N_0^0}{\sqrt{H_0}}X_0^T\sqrt{R_z}[X_{1v_1^1}\
\ ...\ \ X_{1v_1^n}]+2\sqrt{z}[0\ N_0^0\ \ e_{n+2}\ \ ...\ \
e_{2n-1}]\frac{M}{|V_1|^2+1})dV_1,\
dV_1=\frac{(|V_1|^2+1)(|V_0|^2+1)}{2\sqrt{z}}R_1\del
R_0^T(X_1^T\sqrt{R_z}[X_{0v_0^1}\ \ ...\ \ X_{0v_0^n}])^T$, where
$M\subset\mathbf{O}_n(\mathbb{C})$ is determined up to
multiplication on the left with an orthogonal matrix which fixes
$e_1$ (due to the indeterminacy of the rolling in the normal
bundle) by $M^Te_1=\pm
i\frac{(|V_1|^2+1)(X_0^T\sqrt{R_z}[X_{1v_1^1}\ \ ...\ \
X_{1v_1^n}])^T}{2\sqrt{zH_0}}$ (because of the ACPIA requirement
$\frac{4z(e_1^TMdV_1)^2}{(|V_1|^2+1)^2}=-[(dN_0^0)^T(x_z^1-x_0^0)]^2=
-\frac{1}{H_0}(X_0^T\sqrt{R_z}[X_{1v_1^1}\ \ ...\ \
X_{1v_1^n}]dV_1)^2$). The vectors of the facet are obtained by
replacing $\del$ with diagonal constant matrices $\Del$; let
$\Del':=\frac{(|V_1|^2+1)(|V_0|^2+1)}{2\sqrt{z}}R_1\Del
R_0^T(X_1^T\sqrt{R_z}[X_{0v_0^1}\ \ ...\ \ X_{0v_0^n}])^T$.  Thus
we need $e_j^TM\Del'=0,\ j=3,...,n$, so
$M\Del'=e_1^TM\Del'e_1+e_2^TM\Del'e_2$; from $(N_z^1)^Tw=0$ we
obtain $e_2^TM\Del'=\pm ie_1^TM\Del'$; after scaling we obtain
$w=\pm[x_{zv_1^1}^1\ \ ...\ \ x_{zv_1^n}^1](e_1\pm iM^Te_2)$. Now
the condition $w^TAR_z^{-1}w=0$ that $w$ is a ruling on $x_z^1$
becomes $|e_1\pm iM^Te_2|^2=0$ and it is thus satisfied.

The W congruence property follows in a manner similar to the
(I)QWC case.

\section{The Bianchi Theorems on confocal quadrics}

\subsection{The Ivory affinity provides a rigid motion}

\subsection{The second iteration of the tangency configuration}

\subsection{M\"{o}bius configurations}

\subsection{Homographies and confocal quadrics}

\end{document}